\documentclass[10pt, a4paper]{amsart}
\usepackage{amsfonts}
\usepackage{amsthm}
\usepackage{amsmath}

\usepackage{amssymb}
\usepackage{latexsym}

\theoremstyle{plain} \newtheorem{lem}{Lemma}[section]
\theoremstyle{plain} \newtheorem{prop}[lem]{Proposition}
\theoremstyle{plain} 
\theoremstyle{plain} \newtheorem{cor}[lem]{Corollary}

\theoremstyle{definition}

\newtheorem{rem}[lem]{Remark}

\newcommand{\Q}{\mathbb{Q}}
\newcommand{\R}{\mathbb{R}}
\newcommand{\C}{\mathbb{C}}
\newcommand{\F}{\mathbb{F}}
\newcommand{\Z}{\mathbb{Z}}
\renewcommand{\P}{\mathbb{P}}
\newcommand{\Aut}{\mathrm{\mathop{Aut}}}

\newcommand{\kar}{\mathrm{\mathop{char}}}
\newcommand{\GL}{\mathrm{\mathop{GL}}}

\begin{document}

\title{Classification of nilpotent associative algebras of small dimension}
\author{Willem A. de Graaf}
\address{Dipartimento di Matematica\\
Universit\`{a} di Trento\\
Italy}
\date{}

\begin{abstract}
We classify nilpotent associative algebras of dimensions up to $4$ over any
field. This is done by constructing the nilpotent associative algebras
as central extensions of algebras of smaller dimension, analogous to
methods known for nilpotent Lie algebras.
\end{abstract}

\maketitle

\section{Introduction}\label{sec:intro}

The classification of associative algebras is an old and often recurring
problem.
The first investigation into it was perhaps done by Peirce (\cite{peirce}).
Many other publications related to the problem have appeared. Without any
claim of completeness, we mention work by Hazlett (\cite{hazlett16}, nilpotent
algebras of dimension $\leq 4$ over $\C$), Kruse and Price (\cite{kruseprice},
nilpotent associative algebras of dimension $\leq 4$ over any field),
Mazzola (\cite{mazzola79} - associative unitary algebras
of dimension 5 over algebraically closed fields of characteristic not 2,
\cite{mazzola80} - nilpotent commutative associative algebras of dimension
$\leq 5$, over algebraically closed fields of characteristic not 2,3) and
Poonen (\cite{poonen} - nilpotent commutative associative
algebras of dimension $\leq 5$, over algebraically closed fields). Recently,
Eick and Moede (\cite{eickmoede}) have developed a coclass theory for
nilpotent associative algebras, offering a different perspective on their
classification.

The purpose of this paper is, similar to the work of Kruse and Price, to
describe a classification of nilpotent associative algebras of dimensions
up to 4 over any field. We rewrite this classification using a more uniform
method (constructing nilpotent associative algebras as central extensions of
algebras of lower dimension), and provide some small corrections to the
classification in \cite{kruseprice}. Although the method employed here is
quite different, many details have been taken from the proofs of Kruse and
Price. This goes especially for the proofs in Section \ref{sec:noncom2d}; 
examples are the invariant $\sigma$ for solving the isomorphism
problem for the algebras $A_{4,23}^{\alpha,\beta,\gamma}$, the substitution
$u=\omega_2\tau +\omega_1\psi\sigma$, $v=\omega_1\psi-1$ in the proof of
Lemma \ref{lem:A21.8}, the group $H_\alpha$ in Lemma \ref{lem:A21.9}.

On many occasions we use the algorithmic technique of Gr\"obner bases, executed
with help of the computer algebra system {\sc Magma} (\cite{Magma}), to
obtain conditions for isomorphism. This may not entirely satisfy the
theoretically minded reader who wants to verify the results without the
help of a computer. However, as illustrated in Section \ref{sec4}, it is
possible to compute polynomials that make such a verification equivalent to
checking arithmetic identities in polynomial rings. Because they can be a
bit bulky we do not include these polynomials here (except for one instance
in Section \ref{sec3}). 

It may not be entirely obvious what exactly is meant by a classification of
4-dimensional nilpotent associative algebras. Let $V$ be a 4-dimensional
vector space over the field $F$. Let $\mathcal{A}$ be the 64-dimensional space
of bilinear maps $V\times V\to V$, that is, $\mathcal{A}$ is the space of
algebra structures on $V$. Define an action of $\GL(V)$ on $\mathcal{A}$ by
$g\cdot \varphi (v,w) = \varphi(g^{-1}v,g^{-1}w)$ (where $g\in \GL(V)$,
$\varphi\in \mathcal{A}$, $v,w\in V$). Then $\varphi,\psi\in \mathcal{A}$ are
isomorphic if and only if they lie in the same $\GL(V)$-orbit.
Let $\mathcal{N}\subset \mathcal{A}$ consist of
the associative and nilpotent algebra structures. Then $\mathcal{N}$ is
Zariski-closed in $\mathcal{A}$; moreover, polynomials defining $\mathcal{N}$
can explicitly be written down. A classification of 4-dimensional nilpotent
associative algebras over $F$ is a map $\Gamma : S \to \mathcal{N}$,
where $S$ is a set,
such that each $\varphi\in \mathcal{N}$ is isomorphic to exactly one element
of $\Gamma(S)$. One such classification is very easy to give. Indeed, let
$S_0\subset\mathcal{N}$ be a fundamental domain for the action of $\GL(V)$ and
let $\Gamma : S_0\to \mathcal{N}$
be the inclusion map. However, it is clear that such a
classification is not very explicit. Firstly, it is not clear what the
set $S_0$ looks like; for instance, if $F$ is finite then $S_0$
is finite, but we have no direct way to determine its size.
Secondly, this classification
does nothing to help us solve the isomorphism problem: given a $\varphi\in
\mathcal{N}$, to which element of $\Gamma(S_0)$ is it isomorphic?

In Section \ref{sec:list} 
a list will be given of 4-dimensional nilpotent associative algebras. Some
of them depend on parameters running through (a subset of) the ground field.
In each such case a precise condition will be given characterizing the
pairs of isomorphic algebras, obtained relative to different values of the
parameters. Using these conditions it is straightforward to define a set
$S$ as above, that is much more explicit than the set $S_0$.
The proof of the correctness of this list is contained in
Sections \ref{sec3}, \ref{sec4}, \ref{sec5}. In Section \ref{sec:fields}
the set $S$ is made completely explicit for the cases where $F$ is
a finite field and $F=\R$. In the case where $q=|F|$ is finite
its size is determined: it is $5q+20$ for $q$ odd and $5q+17$ for $q$ even.
This is confirmed by experimental data obtained by Moede for $q$ up to
32, using the
{\sf ccalgs} package (\cite{ccalgs}) for {\sf GAP}4 (\cite{gap4}).
Moreover, the proof of the correctness of the list
of Section \ref{sec:list} is constructive: for a given 4-dimensional nilpotent
associative algebra it is possible, by following the steps in the proof, to
find the algebra of the list to which it is isomorphic. In Section
\ref{sec:isompb} this is illustrated in a small example.

{\bf Acknowledgement:} I thank Andrea Caranti for suggesting this problem
to me. Also I am grateful to Heiko Dietrich for a useful email exchange on
the proof of Lemma \ref{lem:charodd2}, and to Tobias Moede for sharing his
computational data with me.

\section{The classification method}\label{sec:method}

The proofs of the main results in this section are simply translations of those for Lie algebras
(cf. \cite{skjelsund}, \cite{gra12}), and are therefore omitted.

Throughout the ground field of the vector spaces and algebras will be denoted $F$. 

\subsection{Central extensions}

Let $A$ be an associative algebra, $V$ a vector space, and $\theta : A\times A
\to V$ a bilinear map. 
Set $A_{\theta} = A\oplus V$. For $a,b\in A$, $v,w\in V$ we define
$(a+v)(b+w)=ab+\theta(a,b)$. Then $A_{\theta}$ is an associative algebra
if and only if 
$$\theta(ab,c) = \theta(a,bc)\text{ for all $a,b,c\in A$}.$$
The bilinear $\theta$ satisfying this are called cocycles. The
set of all cocycles is denoted $Z^2(A,V)$. The algebra $A_{\theta}$
is called a ($\dim V$-dimensional) central extension of $A$ by $V$ (note that 
$A_\theta V=VA_\theta=0$). \par
Let $\nu : A\to V$ be a linear map, and define $\eta(a,b)=\nu(ab)$.
Then $\eta$ is a cocycle, called a coboundary. The set of all coboundaries
is denoted $B^2(A,V)$. Let $\eta$ be a coboundary; then 
$A_{\theta}\cong A_{\theta+\eta}$. Therefore we consider the set 
$H^2(A,V)=Z^2(A,V)/B^2(A,V)$. If we view $V$ as a trivial $A$-bimodule, then
$H^2(A,V)$ is the second Hochschild-cohomology space (cf. \cite{rsp}).
\par
Now let $B$ be an associative algebra. By $C(B)$ we denote the ideal
consisting of all $b\in B$ with $bB=Bb=0$. This is called the {\em annihilator}
of $B$. Suppose that $C(B)$ is 
nonzero, and set $V=C(B)$, and $A=B/C(B)$. Then there is a $\theta\in H^2(A,V)$
such that $B\cong A_\theta$.

We conclude that any  algebra with a nontrivial annihilator
can be obtained
as a central extension of a  algebra of smaller dimension. So in
particular, all nilpotent  algebras can be constructed this way.

When constructing nilpotent algebras as $A_{\theta}=A\oplus V$,
we want to restrict to $\theta$ such that $C(A_{\theta})=V$. If
the annihilator
of $A_{\theta}$ is bigger, then it can be constructed as
a central extension of a different  algebra. (This way we avoid
constructing the same  algebra as central extension of different 
algebras.) Now the radical of a $\theta\in Z^2(A,V)$ is
$$\theta^{\perp} =\{ a\in A\mid \theta(a,b)=\theta(b,a)=0\text{ for all $b\in A$}\}.$$

Then $C(A_{\theta}) = 
(\theta^{\perp}\cap C(A))+V$, proving the following proposition.

\begin{prop}\label{prop:1}
$\theta^{\perp}\cap C(A) =0$ if and only if $C(A_{\theta})=V$.
\end{prop}

Let $e_1,\ldots,e_s$ be a basis of $V$, and $\theta \in Z^2(A,V)$. Then
$$\theta(a,b) = \sum_{i=1}^s \theta_i(a,b) e_i,$$
where $\theta_i \in Z^2(A,F)$. Furthermore, $\theta$ is
a coboundary if and only if all $\theta_i$ are.\par
Let $\phi\in \Aut(A)$. For $\eta\in Z^2(A,V)$ define $\phi\eta(a,b) = 
\eta(\phi(a),\phi(b))$. Then $\phi\eta \in Z^2(A,V)$. So $\Aut(A)$ acts on 
$Z^2(A,V)$. Also, $\eta\in B^2(A,V)$ if and only if $\phi\eta\in B^2(A,V)$
so $\Aut(A)$ acts on $H^2(A,V)$. 

\begin{prop}\label{prop:2}
Let $\theta(a,b) = \sum_{i=1}^s \theta_i(a,b)e_i$ and 
$\eta(a,b)=\sum_{i=1}^s \eta_i(a,b) e_i$. Suppose that $\theta^{\perp}
\cap C(A)= \eta^{\perp}\cap C(A)=0$. Then
$A_{\theta}\cong A_{\eta}$ if and only if there is a $\phi\in\Aut(A)$
such that the $\phi\eta_i$ span the same subspace of $H^2(A,F)$
as the $\theta_i$.
\end{prop}

Let $A=I_1\oplus I_2$ be the direct sum of two ideals. Suppose that 
$I_2$ is contained in the annihilator of $A$. 
Then $I_2$ is called a {\em central component} of $A$. 

\begin{prop}\label{prop:3}
Let $\theta$ be such that $\theta^{\perp}\cap C(A) =0$. Then
$A_{\theta}$ has no central components if and only if $\theta_1,
\ldots,\theta_s$ are linearly independent.
\end{prop}

Based on Propositions \ref{prop:1}, \ref{prop:2}, \ref{prop:3}
we formulate a procedure that takes as input a nilpotent
algebra $A$ of dimension $n-s$. It outputs all nilpotent algebras $B$
of dimension $n$ such that $B/C(B)\cong A$, and $B$ has no central components.
For this we need some more terminology. Let $\Omega$ be an $s$-dimensional 
subspace of $H^2(A,F)$ spanned by $\theta_1,\ldots,\theta_s$. Let $V$ be 
an $s$-dimensional vector space spanned by $e_1,\ldots,e_s$. Then we define
$\theta\in H^2(A,V)$ by $\theta(a,b) = \sum_i \theta_i(a,b)e_i$. We call $\theta$
the cocycle corresponding to $\Omega$ (or more precisely, to the chosen basis of $\Omega$).
Furthermore, we say that $\Omega$ is 
{\em useful} if $\theta^\perp \cap C(A)=0$. Note that $\theta^\perp$ is the 
intersection of the $\theta_i^\perp$.

Now the procedure runs as follows:

\begin{enumerate}
\item Determine $Z^2(A,F)$, $B^2(A,F)$ and $H^2(A,F)$.
\item Determine the orbits of $\Aut(A)$ on the set of useful $s$-dimensional
subspaces of $H^2(A,F)$.
\item For each orbit let $\theta$ be the cocycle corresponding to a representative
of it, and construct $A_\theta$. 
\end{enumerate}

Of course the hard part is Step 2. Note that $\Aut(A)$ is an
algebraic group. This means that whether two useful subspaces lie in the
same $\Aut(A)$-orbit is equivalent to the existence of a solution over $F$
of a set of polynomial equations. On some occasions we cannot decide solvability by hand.
Then we use the technique of Gr\"obner bases (cf. \cite{clo}). This is an
algorithmic procedure to compute an equivalent set of polynomial equations
that is sometimes easier to solve. On all occasions where we use this the
equations have coefficients in $\Z$. For the Gr\"obner basis calculation
we take the ground field to be $\Q$. A priori this yields results that are
only valid over fields of characteristic 0. However, the {\sc Magma} computational
algebra system (\cite{Magma}) has the facility to compute the coefficients of 
an element of the Gr\"obner basis relative to the input basis. We use this in
order to derive conclusions valid in all characteristics. This will be illustrated in
more detail in Section \ref{sec4}.

The procedure only gives those algebras without central components. So we have
to add the algebras obtained by taking the direct sum of a smaller-dimensional
algebra with a nil-algebra (that has trivial multiplication).

\subsection{Notation and terminology}

The base field of all algebras will be denoted $F$. Furthermore, $F^*$ is the
set of nonzero elements of $F$.

Let $A$ be an associative algebra with basis elements $a_1,\ldots,a_n$.
Then by $\Delta_{a_i,a_j}$ we denote the bilinear map $A\times A\to F$
with $\Delta_{a_i,a_j} (a_k,a_l)  = 1$ if $i=k$ and $j=l$, and otherwise
it takes the value $0$.

Throughout the basis elements of the algebras will be denoted by the letters
$a,b,\ldots$. We specify an algebra by an expression in angled brackets. First
we list the basis elements, and then the nonzero products among the basis
elements. For example:
$$A = \langle a,b,c \mid a^2=b^2=c \rangle$$
specifies the algebra with basis $a,b,c$, with $a^2=b^2=c$ and the other
products among the basis elements are zero.

Let $A$ be a nilpotent associative algebra, $M=H^2(A,F)$, $G=\Aut(A)$.
The main problem that we will be dealing with is to list the orbits of $G$ on
the set of $s$-dimensional usable subspaces of $M$. (In fact, we will always
have $s=1$ or $s=2$.) These subspaces correspond to points in the Grassmannian
$\mathrm{Gr}(M,s)$ of $s$-dimensional subspaces of $M$. We say that two such
points are \emph{conjugate} if they lie in the same $G$-orbit.

Often we will be dealing with a Grassmannian of 1-dimensional
subspaces. In that situation we say that $\theta_1,\theta_2\in M$ are conjugate
if there is a $\phi\in G$ and a $\lambda \in F^*$ such that $\phi\theta_1 =
\lambda\theta_2$.

\section{The list of nilpotent associative algebras of dimension 4}
\label{sec:list}

In this section we give the list of 4-dimensional nilpotent associative
algebras. The list consists of single algebras and of parametrized series
of algebras. The latter are followed by restrictions on the parameters, on the
ground field, and by a precise description of the isomorphisms that exist
between different members of the series. The algebras $A_{4,k}$, $A_{4,l}$ for
$k\neq l$ are not isomorphic. 

\begin{itemize}
\item $A_{4,1}=\langle a,b,c,d \mid \rangle$.
\item $A_{4,2}=\langle a,b,c,d \mid  a^2 = b\rangle$.
\item $A_{4,3}^\delta=\langle a,b,c,d\mid a^2 = c,~ b^2 = \delta c\rangle$, 
  $\delta \neq 0$, $A_{4,3}^\delta\cong A_{4,3}^\epsilon$ if and only if there is a
  $\nu\in F^*$ with $\delta = \nu^2 \epsilon$.   
\item $A_{4,4}^\delta=\langle a,b,c,d \mid a^2 = c,~ b^2 = \delta c, ~
  ab=c\rangle$,
\item $A_{4,5}=\langle a,b,c,d \mid ab=c, ~ ba = -c \rangle$.
\item $A_{4,6}=\langle a,b,c,d\mid a^2=b, ab=ba = c\rangle$.
\item $A_{4,7} = \langle a,b,c,d \mid a^2=bc=-cb = d\rangle$, $\kar(F)\neq 2$.
\item $A_{4,8}^{\alpha,\beta} = \langle a,b,c,d \mid a^2=d,~ b^2=\alpha d,~ c^2 =
  \beta d \rangle$, where $\alpha,\beta\neq 0$. We have that
  $A_{4,8}^{\alpha,\beta}\cong A_{4,8}^{\gamma,\delta}$ if and only if if and only if
  the quadratic forms
  $\alpha x^2+\beta y^2 +\alpha \beta z^2$, $\gamma x^2 +\delta y^2+\gamma\delta
  z^2$ are equivalent.
\item $A_{4,9}^{\alpha,\beta} = \langle a,b,c,d \mid a^2=\alpha d,~ b^2=d,~ bc=d,~
  c^2 = \beta d \rangle$, where  $\alpha\neq 0$,  $A_{4,9}^{\alpha,\beta}
  \cong A_{4,9}^{\gamma,\delta}$
($\alpha,\gamma\neq 0$) if and only if $\beta=\delta$ and there are $s,t\in F$
with $t^2-st+\delta s^2 = \tfrac{\alpha}{\gamma}$.
\item $A_{4,10}^{\alpha} = \langle a,b,c,d \mid a^2= d,~ ab=d,~ b^2=\alpha d,~
  bc=d,~cb = d \rangle$; if $\kar(F)\neq 2$ then $A_{4,10}^\alpha
\cong A_{4,10}^0$; if $\kar(F)=2$ then $A_{4,10}^{\alpha} \cong
A_{4,10}^\beta$ if and only if there is a $T\in F$ with $T^2+T+\alpha+\beta=0$.
\item  $A_{4,11} = \langle a,b,c,d \mid a^2= d,~ ab=d,~ cb = d, ~ c^2=-d
  \rangle$.
\item $A_{4,12} = \langle a,b,c,d \mid a^2=b, ab=ba=d, ac=d, c^2=d \rangle$.
\item $A_{4,13} = \langle a,b,c,d \mid a^2=b, ab=ba=d, c^2=d \rangle$.
\item $A_{4,14} = \langle a,b,c,d \mid a^2=b, ab=ba=d, ac=d \rangle$.
\item $A_{4,15} = \langle a,b,c,d \mid a^2=b, ab=ba=c, b^2=ac=ca=d \rangle$.
\item $A_{4,16} = \langle a,b,c,d \mid a^2=c, ba=d \rangle$.
\item $A_{4,17}^\delta = \langle a,b,c,d \mid a^2=c, ab=ba=d, b^2=\delta c +
  d\rangle$, $\kar(F)=2$, $A_{4,17}^\delta \cong A_{4,17}^\epsilon$ if and only
  if there is a $T\in F$ with $T^2+T+\delta+\epsilon=0$.
\item $A_{4,18}^\delta = \langle a,b,c,d \mid a^2=c, ab=d, ba=-d, b^2=\delta c
  \rangle$. If $\kar(F)\neq 2$ then $A_{4,18}^\delta \cong A_{4,18}^\epsilon$ if
  and only if there is a $\nu\in F^*$ with $\epsilon = \nu^2 \delta$. If
  $\kar(F)=2$ $A_{4,18}^\delta \cong A_{4,18}^\epsilon$ if and only if there are
  $u,v,x,y\in F$ with $uy+vx\neq 0$, $u^2+v^2\delta \neq 0$ and
  $\epsilon = \tfrac{x^2+y^2\delta}{u^2+v^2\delta}$.
\item $A_{4,19}^\delta = \langle a,b,c,d \mid a^2=c, ab=d, ba=c+d, b^2=\delta
  c\rangle$; $\kar(F)=2$ and $A_{4,19}^\delta\cong A_{4,19}^\epsilon$ if and only if
  there is a $T\in F$ with $T^2+T+\delta+\epsilon=0$.
\item $A_{4,20} = \langle a,b,c,d \mid a^2=c, ab=d, ba=c \rangle$.
\item $A_{4,21}^\delta = \langle a,b,c,d \mid a^2=c, ab=d, ba=\delta d \rangle$,
  $\delta\neq -1$.
\item $A_{4,22} = \langle a,b,c,d \mid a^2=c, ab=d, ba=c+d, b^2=c \rangle$,
  $F=\F_3$.
\item $A_{4,23}^\delta = \langle a,b,c,d \mid a^2=c, ab=ba=d, b^2=-\delta
  c\rangle$, $\delta\neq 0$, $\kar(F)\neq 2$ and
  $A_{4,23}^\delta \cong A_{4,23}^\epsilon$ if
  and only if there is a $\nu\in F^*$ with $\epsilon = \nu^2 \delta$.
\item $A_{4,24} = \langle a,b,c,d \mid a^2=c, ab=d, ba=-c, b^2=c \rangle$,
  $\kar(F)\neq 2$.
\item $A_{4,25}^{\alpha,\beta,\gamma} = \langle a,b,c,d \mid a^2=c, ab=d,
  ba=-\beta c+\alpha d, b^2= -\gamma c \rangle$, $\kar(F)\neq 2$,
  $\alpha\neq \pm 1$,  $\gamma\neq 0$, $\sigma^2\neq -\gamma$, where $\sigma =
  \tfrac{\beta}{1-\alpha}$.  Furthermore,
  $A_{4,25}^{\alpha,\beta,\gamma}\cong A_{4,25}^{\alpha',\beta',\gamma'}$ if and only if
  there are $\varphi,\nu\in F^*$ with
  $(\sigma')^2+\gamma' = \varphi^2(\sigma^2+\gamma)$ (where $\sigma' =
  \tfrac{\beta'}{1-\alpha'}$), $\gamma = \nu^2\gamma'$
  and, setting $\psi = \nu\varphi$, letting $\omega_1=\pm 1$ be such that
  $-\omega_1\neq \tfrac{1+\alpha}{1-\alpha}\psi$ and defining
  $$\lambda = \frac{\psi(1+\alpha)-\omega_1(1-\alpha)}{
  \psi(1+\alpha) + \omega_1(1-\alpha)}\text{ and }
  \mu = \sigma'(1-\lambda),$$
  we have $\alpha'=\lambda$, $\beta'=\mu$ or $\alpha' = \lambda^{-1}$, $\beta'=
  -\mu \lambda^{-1}$.
\item $A_{4,26}^{\alpha,\beta,\gamma} = \langle a,b,c,d \mid a^2=c, ab=d,
  ba=\beta c+\alpha d, b^2= \gamma c \rangle$, $\kar(F)= 2$,
  $\alpha\neq 1$,  $\gamma\neq 0$. Furthermore,
  $A_{4,26}^{\alpha,\beta,\gamma}\cong A_{4,26}^{\alpha',\beta',\gamma'}$ if and only if
  there is a $\nu\in F^*$ with $\gamma = \nu^2\gamma'$ and $\sigma =
  \nu \sigma'$ (where $\sigma = \tfrac{\beta}{1+\alpha}$, 
  $\sigma' = \tfrac{\beta'}{1+\alpha'}$), and a $h\in H_{\alpha,\beta,\gamma}$
  with $\tfrac{1}{1+\alpha'} = \tfrac{1}{1+\alpha} +h$, where
  $$H_{\alpha,\beta,\gamma} = \{ \frac{\sigma uv + \gamma v^2}{u^2+\gamma v^2} \mid
  u,v\in F \text{ and } u^2+\gamma v^2\neq 0 \},$$
  which is an additive subgroup of $F$.
\end{itemize}

\begin{rem}
The algebras $A_{4,k}$, $1\leq k\leq 6$ and $A_{4,17}^0$, $A_{4,23}^{-1}$ are
decomposable (i.e., they are direct sums of ideals). The others are not.  
\end{rem}

\begin{rem}
  As remarked in Section \ref{sec:intro}, many details of the proof of the
  correctness of this list have been taken from \cite{kruseprice}. There
  are, however, also some differences in the final result. It appears that
  in \cite{kruseprice} it is stated that $A_{4,8}^{\alpha,\beta}\cong
  A_{4,8}^{\gamma,\delta}$ if and only if the quadratic forms $x^2+\alpha y^2
  +\beta z^2$, $x^2+\gamma y^2+\delta z^2$ are equivalent, which is different
  from the condition obtained here. Furthermore, the
  algebras $A_{4,8}^{\alpha,\beta}$ such that $x^2+\alpha y^2+\beta z^2$ has a
  nontrivial zero, are omitted. I can only explain that by the fact that
  in \cite{kruseprice} only the indecomposable algebras are classified.
  However, the algebras $A_{4,8}^{\alpha,\beta}$ are not decomposable, regardless
  of the existence of a zero of $x^2+\alpha y^2+\beta z^2$.
  In \cite{kruseprice} it is claimed that $A_{4,9}^{\alpha,\beta}\cong
  A_{4,9}^{\gamma,\delta}$ if and only if $\beta=\delta$ and $\alpha = \nu^2 \gamma$
  for some $\nu\in F^*$. This is sufficient but not necessary. Finally, in
  \cite{kruseprice} the classification is presented using sets that precisely
  parametrize the non-isomorphic algebras. For example, we have $A_{4,3}^\delta$
  for $\delta \in F^*/(F^*)^2$. Here we have stated the conditions under
  which two algebras are isomorphic, as this helps in solving the isomorphism
  problem and also immediately characterizes the parameter sets (like
  $F^*/(F^*)^2$. For this we
  have taken the analysis of the isomorphism of the algebras
  $A_{4,23}^{\alpha,\beta,\gamma}$ one step further than in \cite{kruseprice}.
\end{rem}  

\section{Dimensions 1 and 2}\label{sec3}

There is only one nilpotent algebra of dimension 1: it is spanned
by $a$, and $a^2=0$. We denote it by $A_{1,1}$.

Now $H^2(A_{1,1},F)$ is spanned by $\Delta_{a,a}$. So we get two nilpotent 
algebras of dimension $2$, corresponding to $\theta=0$ and $\theta=
\Delta_{a,a}$ respectively. They are $A_{2,1}$, which is spanned by $a,b$,
and all products are zero, and 
$$A_{2,2} = \langle a,b \mid a^2 = b\rangle .$$

\section{Dimension 3}\label{sec4}

In this section we classify nilpotent associative algebras of dimension 
3, over any field.

First we get the algebras that are the direct sum of an algebra of dimension
$2$ and a $1$-dimensional ideal, isomorphic to $A_{1,1}$, spanned by $c$.
We denote them $A_{3,1}$ (all products zero), and 
$$A_{3,2}=\langle a,b,c \mid  a^2 = b\rangle.$$

There are no $2$-dimensional central extensions of $A_{1,1}$.
So we consider $1$-dimensional central extensions of $A_{2,1}$. Here 
$H^2(A_{2,1},F)$ consists of $\theta = \alpha\Delta_{a,a}+\beta\Delta_{a,b}
+\gamma \Delta_{b,a} + \delta \Delta_{b,b}$. The automorphism group
consists of all 
$$\phi=\begin{pmatrix} u & x \\ v & y \end{pmatrix}, \text{ with } 
uy-vx \neq 0.$$
Write $\phi\theta = \alpha' \Delta_{a,a} + \cdots +\delta' \Delta_{b,b}$.
Then
\begin{align*}
\alpha' &= u^2 \alpha +uv\beta + uv\gamma +v^2\delta\\
\beta' &= ux\alpha + uy\beta + vx\gamma + vy\delta\\
\gamma' &= ux\alpha +vx\beta +uy\gamma +vy\delta\\
\delta' &= x^2\alpha + xy\beta + xy\gamma +y^2\delta.
\end{align*}

We distinguish a few cases.

{\em Case 1:} suppose that there are $f\in A_{2,1}$ with $\theta(f,f)\neq 0$.
Then we may assume that $\alpha\neq 0$, and we can divide to get $\alpha =1$.
Choose $u=y=1$, $v=0$, $x=-\gamma$ to get $\gamma'=0$ and $\alpha'=1$.
So we may assume that $\alpha=1$ and $\gamma=0$. Choose $x=v=0$, $u=1$;
this leads to $\alpha'=1$, $\gamma'=0$, and $\beta' = y\beta$. We can
still freely choose $y\neq 0$. So we are left with two cases: $\beta = 0,1$.

{\em Case 1a.}
If $\beta = 0$, then we get the cocycles $\theta^1_\delta = \Delta_{a,a}
+\delta\Delta_{b,b}$. If we choose $x=v=0$, $u=1$, then $\alpha'=1$, $\gamma'=
\beta'=0$ and $\delta' = y^2 \delta$. So we see that $\theta^1_\delta$ and
$\theta^1_{y^2 \delta}$ are conjugate for any $y\neq 0$. 
In order to show the converse let $\phi$ be as above. Then 
$\phi\theta^1_\delta = \lambda \theta^1_\epsilon$ (for some $\lambda\in
F^*$) amounts to the following polynomial equations
\begin{align*}
f_1 &:= u^2 +v^2\delta-\lambda =0,\\
f_2 &:= ux +vy\delta =0,\\
f_3 &:= x^2+y^2\delta-\lambda\epsilon=0.
\end{align*}

To these we add 
$$ f_4:= D(uy-vx)-1 =0,$$
which ensures that $\det\phi\neq 0$. 

Now a reduced Gr\"obner basis of the ideal generated by $f_1,\ldots ,f_4$
contains the polynomials $u^2\epsilon-y^2\delta$, and $v^2\delta\epsilon-x^2$.
Using {\sc Magma} it is not only possible to compute this Gr\"obner basis,
but also to write its elements in terms of the $f_i$. In this case we have
\begin{align*}
  u^2\epsilon -y^2\delta &= (Dvx\epsilon+\epsilon)f_1+(Dxy\beta-Duv\epsilon)f_2-
  (Dvx\beta+1)f_3+(v^2\delta\epsilon-x^2)f_4\\
  v^2\delta\epsilon-x^2 &= -Dvx\epsilon f_1+(Duv\epsilon-Dxy)f_2+Dvx f_3
  -(v^2\delta\epsilon-x^2)f_4.
\end{align*}

We see that the coefficients that appear all lie in $\Z$; so these equations are valid
over any field $F$.
Hence if there is a $\phi\in \Aut(A_{2,1})$ with $\phi\theta^1_\delta =
\lambda \theta^1_\epsilon$,
then there are $u,v,x,y\in F$ with $u^2\epsilon-y^2\delta = v^2\delta\epsilon-
x^2=0$ and $uy-vx\neq 0$. This implies that there exists $y\in F^*$
with $\delta = y^2 \epsilon$. The conclusion is that 
$\theta^1_\delta$ and $\theta^1_\epsilon$ are conjugate 
if and only if there is a $y\in F^*$
with $\delta = y^2\epsilon$.

{\em Case 1b.}
If $\beta=1$, then we get the cocycles $\theta^2_\delta = \Delta_{a,a}
+\Delta_{a,b} + \delta\Delta_{b,b}$. Note that these cannot be
$\Aut(A_{2,1})$-conjugate to a 
$\theta^1_\delta$ as the latter is symmetric. Also here we use a Gr\"obner basis 
calculation, of which we do not give all the details.
In this case when we write the polynomial equations that are equivalent to
$\phi\theta^2_\delta = \lambda \theta^2_\epsilon$ and compute a 
Gr\"obner basis, then we find that it contains $\delta-\epsilon$. Also, writing
$\delta-\epsilon$ in terms of the initial polynomials (as above) we conclude
that this is valid over all fields. So, in this case 
$\theta^2_\delta$ and $\theta^2_\epsilon$ are conjugate 
if and only if $\delta=\epsilon$.

{\em Case 2:} $\theta(f,f)=0$ for all $f\in A_{2,1}$. In that case, $\alpha=
\delta=0$ and $\beta=-\gamma$. So, after dividing we may assume $\beta=1$,
$\gamma=-1$, and we
get $\theta^3 = \Delta_{a,b} - \Delta_{b,a}$. We have 
that $\phi\theta^3$ is a multiple of $\theta^3$. Hence it is not conjugate
to any of the previous cocycles.

So we get the nonzero cocycles $\theta^1_\delta$, $\theta^2_\delta$,
and $\theta^3$. For the first we need
$\delta\neq 0$, otherwise $b$ lies in the radical. This leads to the algebras: 
$$A_{3,3}^\delta=\langle a,b,c\mid a^2 = c,~ b^2 = \delta c\rangle,~~~
\delta \neq 0,$$
$$A_{3,4}^\delta=\langle a,b,c \mid a^2 = c,~ b^2 = \delta c, ~ ab=c\rangle,$$
$$ A_{3,5}=\langle a,b,c \mid ab=c, ~ ba = -c \rangle.$$
From the above discussion it follows that $A_{3,3}^\delta$ is isomorphic to
$A_{3,3}^\epsilon$ if and only if there is an $y\in F^*$ with $\delta = y^2 
\epsilon$. 

Next we consider $1$-dimensional central extensions of $A_{2,2}$.
Here we get that $Z^2(A_{2,2},F)$ is spanned by $\Delta_{a,a}$ and
$\Delta_{a,b}+\Delta_{b,a}$. Moreover, $B^2(A_{2,2},F)$ is spanned
by $\Delta_{a,a}$. So we get only one cocycle $\theta =
\Delta_{a,b}+\Delta_{b,a}$, yielding the algebra
$$A_{3,6}=\langle a,b,c\mid a^2=b, ab=ba = c\rangle.$$

Concluding, we have the following nilpotent 3-dimensional algebras:
$A_{3,1}$, $A_{3,2}$, $A_{3,3}^\delta$, where $\delta\in F^*/(F^*)^2$,
$A_{3,4}^\delta$,
where $\delta\in F$, $A_{3,5}$, $A_{3,6}$. So over an infinite field there is an
infinite number of them, wheras over $\F_q$ there are $q+6$ for $q$ odd,
and $q+5$ for $q$ even.

\begin{rem}
  By inspection it is seen that we have obtained the same classification
  as in \cite{kruseprice}, Theorem 2.3.6.
\end{rem}
 
\section{Nilpotent algebras of dimension 4}\label{sec5}

First we get the algebras that are the direct sum of a $3$-dimensional
algebra, and $A_{1,1}$. This way we get the algebras $A_{4,i}$, $1\leq i\leq 6$.

Next we consider the $1$-dimensional central extensions of the algebras
$A_{3,i}$, $1\leq i\leq 6$. Staightforward calculations show that
a $\theta\in Z^2(A_{3,i},F)$, for $i=3,4,5$, alwas has $c\in \theta^\perp$.
So those algebras do not yield anything. For each remaining case we have
a subsection.

Finally, $A_{2,2}$ does not have 2-dimensional central extensions, so we
are left with determining the 2-dimensional central extensions of $A_{2,1}$,
which is done in Section \ref{sec:noncom2d}.

\subsection{1-dimensional central extensions of $A_{3,1}$}\label{sec:A31}

Let $B=(e_1,e_2,e_3)$ be an ordered basis of $A_{3,1}$.  
Then $H^2(A_{3,1},F)$ consists of all $\theta=\sum_{i,j=1}^3 \gamma_{i,j}
\Delta_{i,j}$, where $\Delta_{ij,}=\Delta_{e_i,e_j}$. We let $[\theta]_B$ denote
the $3\times 3$-matrix $(\gamma_{i,j})$. To ease notation a bit, 
on many occasions we will just identify $\theta$ with $[\theta]_B$. 

We have that $\theta_1,\theta_2\in
H^2(A_{3,1},F)$ are conjugate if and only if
there is a basis $B'$ of $A_{3,1}$ and a
nonzero $\lambda\in F$ with $[\theta_1]_B = \lambda [\theta_2]_{B'}$.
This is equivalent to the existence of a nonsingular $3\times 3$-matrix $M$ with
$M[\theta_1]_BM^T = \lambda [\theta_2]_{B}$.

Let $\theta\in H^2(A_{3,1},F)$.
We distinguish a few cases. {\em Case 1:} $\theta(a,a)=0$ for all $a\in
A_{3,1}$. This means that $\theta$ is an {\em alternate} bilinear form.
By \cite{jac2}, Chapter V, Theorem 7, there is a basis $B$ of $A_{3,1}$
such that $[\theta]_B$ is block diagonal with blocks
$\left(\begin{smallmatrix} 0 & 1\\ -1 & 0\end{smallmatrix}\right)$, or 0.
Hence $\theta$ has a nonzero radical. Therefore the space spanned by $\theta$
is not useful.  

{\em Case 2:} there are $a\in A_{3,1}$ with $\theta(a,a)\neq 0$. Then there is
a basis $B=(e_1,e_2,e_3)$ of $A_{3,1}$ with $\theta(e_1,e_1)\neq 0$. After
dividing, we may assume that $\theta(e_1,e_1)=1$. As above, let $\gamma_{i,j} =
\theta(e_i,e_j)$. Set $e_1'=e_1$, $e_2'=e_2-\gamma_{2,1}e_1$, $e_3'=e_3-
\gamma_{3,1}e_1$. Then
$\theta(e_1',e_1')=1$ and $\theta(e_2',e_1') =\theta(e_3',e_1')=0$. So we may
assume that $\gamma_{2,1}=\gamma_{3,1}=0$.

{\em Case 2a:} $\gamma_{1,2}=\gamma_{1,3}=0$. Let $U$ be the subspace of $A_{3,1}$
spannned by $e_2,e_3$. {\em Case 2aa:} $\theta(u,u)=0$ for all $u\in A_{3,1}$.
Then $\gamma_{2,2}=\gamma_{3,3}=0$ and $\gamma_{2,3}=-\gamma_{3,2}=\alpha$.
We may assume that $\alpha\neq 0$, as otherwise $\theta$ has nonzero radical.
Set $e_1'=\alpha e_1$, $e_2'=\alpha e_2$, $e_3'= e_3$. The matrix of $\theta$
with respect to this basis is $\alpha^2$ times
$$\theta^1=\left(\begin{smallmatrix} 1 & 0 & 0 \\ 0 & 0 & 1 \\ 0 & -1 & 0
\end{smallmatrix}\right),$$
yielding the algebra $A_{4,7}$.

{\em Case 2ab:} there are $u\in U$ with $\theta(u,u)\neq 0$. Then we may assume
that $\gamma_{2,2}\neq 0$. By setting $e_1'=e_1$, $e_2'=e_2$,
$e_3'=\gamma_{3,2}e_2-\gamma_{2,2} e_3$ we see that we may assume that
$\gamma_{3,2}=0$. If $\gamma_{2,3}=0$ as well then we have the cocycles
$$\theta^2_{\alpha,\beta}=\left(\begin{smallmatrix} 1 & 0 & 0 \\ 0 & \alpha & 0
  \\ 0 & 0 & \beta
\end{smallmatrix}\right), \text{ where } \alpha,\beta\neq 0.$$
giving the algebras $A_{4,8}^{\alpha,\beta}$.

If $\gamma_{2,3}\neq 0$ then we set $e_1'=e_1$,
$e_2'= \gamma_{2,3}^{-1} \gamma_{2,2} e_2$, $e_3'=e_3$, showing that we may suppose
that $\gamma_{2,3}=\gamma_{2,2}$. After dividing by $\gamma_{2,2}$ we get
$$\theta^3_{\alpha,\beta}=\left(\begin{smallmatrix} \alpha & 0 & 0 \\
  0 & 1 & 1  \\ 0 & 0 & \beta
\end{smallmatrix}\right), \text{ where } \alpha\neq 0,$$
which gives the algebras $A_{4,9}^{\alpha,\beta}$.

We now consider the conjugacy relations between the cocycles we have
obtained thus far.

If the characteristic is not 2, then $\theta^1$ is not conjugate to
$\theta^2_{\alpha,\beta}$ (as the latter is symmetric), or $\theta^3_{\alpha,\beta}$
(this is seen by a Gr\"obner basis computation). However, if the characteristic
is 2, then by setting
$e_1'=e_1+e_2$, $e_2'=e_1+e_3$, $e_3'=e_1+e_2+e_3$ it is seen that $\theta^1$ is
conjugate to $\theta^2_{1,1}$.

Since $\theta_{\alpha,\beta}^2$ is symmetric, it is not conjugate to
$\theta^3_{\gamma,\delta}$. However, among the $\theta_{\alpha,\beta}^2$ there can
be conjugate pairs, as explained by the following lemma. For the terminology
and notation relative to quadratic forms and quaternion algebras we refer to
\cite{omeara}.

\begin{lem}\label{lem:quad}
  Let $\alpha,\beta,\gamma,\delta \in F^*$. Then $\theta_{\alpha,\beta}^2$ and
  $\theta_{\gamma,\delta}^2$ are conjugate if and only if the quadratic forms
  $\alpha x^2+\beta y^2 +\alpha \beta z^2$, $\gamma x^2 +\delta y^2+\gamma\delta
  z^2$ are equivalent. If the characteristic is not 2 then this holds
  if and only if the quaternion algebras
  $\left( \tfrac{-\alpha,-\beta}{F} \right)$,
  $\left( \tfrac{-\gamma,-\delta}{F} \right)$ are isomorphic.
\end{lem}

\begin{proof}
  We start by showing the first equivalence. Write $X_{\alpha,\beta} =
  \tfrac{1}{\alpha\beta} \theta^2_{\alpha,\beta}$. Suppose that there is a
  nonsingular $3\times 3$-matrix $M$ and $\lambda\in F^*$ with
  $M\theta^2_{\alpha,\beta} M^T = \lambda \theta^2_{\gamma,\delta}$. By taking
  determinants it follows that $\lambda \tfrac{\gamma\delta}{\alpha\beta} =
  \nu^2$ for some $\nu  \in F^*$. Set $N = \tfrac{1}{\nu}M$, then
  $N X_{\alpha,\beta} N^T = X_{\gamma,\delta}$, implying that $\tfrac{1}{\alpha} x^2
  +\tfrac{1}{\beta} y^2 +\tfrac{1}{\alpha\beta} z^2$, $\tfrac{1}{\gamma} x^2
  +\tfrac{1}{\delta} y^2 +\tfrac{1}{\gamma\delta} z^2$ are equivalent.
  Obviously, these two quadratic forms are equivalent to the ones given in the
  lemma.

  Conversely, if  $\alpha x^2+\beta y^2 +\alpha \beta z^2$,
  $\gamma x^2 +\delta y^2+\gamma\delta z^2$ are equivalent then there is
  a nonsingular $3\times 3$-matrix $M$ with $MX_{\alpha,\beta} M^T =
  X_{\gamma,\delta}$. But that implies that $M\theta^2_{\alpha,\beta} M^T =
  \tfrac{\alpha\beta}{\gamma\delta} \theta^2_{\gamma,\delta}$.

  The second equivalence follows from \cite{omeara}, 57:8.
\end{proof}  

We have that $\theta^3_{\alpha,\beta}$ is conjugate to $\theta^3_{\gamma,\delta}$ 
($\alpha,\gamma\neq 0$) if and only if $\beta=\delta$ and there are $s,t\in F$
with $t^2-st+\delta s^2 = \tfrac{\alpha}{\gamma}$. The necessity of this
condition is readily established by a Gr\"obner basis computation. Conversely,
let $B=(e_1,e_2,e_3)$ be a basis of $A_{3,1}$ such that
$$[\theta_{\alpha,\beta}]_B = \left( \begin{smallmatrix} \alpha & 0 & 0 \\
  0 & 1 & 1 \\ 0 & 0 & \beta \end{smallmatrix} \right).$$
Suppose that $\beta=\delta$ and let $s,t\in F$ be given satisfying the above
condition. Set $e_1' = e_1$, $e_2'=(t-s)e_2+s e_3$, $e_3'=-s\delta e_2 + te_3$.
With $B'=(e_1',e_2',e_3')$ we have that $[\theta_{\alpha,\beta}]_{B'} =
\tfrac{\alpha}{\gamma} [\theta_{\gamma,\delta}]_B$.
If the characteristic of $F$ is not 2, then by viewing $t^2-st+\delta s^2 =
\tfrac{\alpha}{\gamma}$ as an equation in $t$, and by considering its
discriminant, one sees that the existence of $s,t$ satisfying this equation,
is equivalent to the existence of $x,y\in F$ with $x^2+(4\delta-1)y^2 = 4
\tfrac{\alpha}{\gamma}$.

{\em Case 2b:} at least one of $\gamma_{1,2},\gamma_{1,3}$ is nonzero. After
possibly interchanging $e_2,e_3$ we may assume that $\gamma_{1,2}\neq 0$.
By setting $e_1'=e_1$, $e_2'= \tfrac{1}{\gamma_{1,2}}e_2$, $e_3'=e_3$ it follows
that we may assume that $\gamma_{1,2}=1$. Then by setting
$e_1'=e_1$, $e_2'=e_2$, $e_3'=e_3-\gamma_{1,3}e_2$ we see that we may assume
that $\gamma_{1,3}=0$.

{\em Case 2ba:} $\gamma_{2,3}= \gamma_{3,2}$. If $\gamma_{2,3}=0$ then set
$e_1'=e_3$, $e_2'=e_1$, $e_3'=e_2$, and we are back in Case 2a. Furthermore,
if $\gamma_{3,3}\neq 0$ then we set $e_1'=e_3$, $e_2'=e_1+\gamma_{3,3}e_2-
\gamma_{2,3} e_3$, $e_3'=\gamma_{3,3}e_2-\gamma_{2,3}e_3$, showing that again we
are back in Case 2a. So we may assume that $\gamma_{2,3}\neq 0$,
$\gamma_{3,3}=0$. Then we set $e_1'=e_1$, $e_2'=e_2$, $e_3'=\tfrac{1}{\gamma_{2,3}}
e_3$, showing that we may assume that $\gamma_{2,3}=1$. We obtain the cocycles
$$\theta^4_{\alpha}=\left(\begin{smallmatrix} 1 & 1 & 0 \\ 0 & \alpha & 1
  \\ 0 & 1 & 0 \end{smallmatrix}\right),$$
yielding the algebras $A_{4,10}^{\alpha}$.
By Gr\"obner basis computations it is seen that $\theta_\alpha^4$ is not
conjugate to $\theta^i$, $i=1,2,3$. 
If the characteristic of $F$ is not 2, then by setting $e_1'=e_1$,
$e_2'=e_2-\tfrac{1}{2}\alpha e_3$, $e_3'=e_3$, we see that $\theta^4_\alpha$ is
conjugate to $\theta^4_0$.
If the characteristic is 2, then $\theta^4_{\alpha}$,  
$\theta^4_\beta$ are conjugate if and only if there is a
$T\in F$ with $T^2+T+\alpha+\beta=0$.
The necessity of this condition is established by a Gr\"obner basis computation.
Conversely, if such a $T\in F$ exists, then set $e_1'=e_1+Te_3$, $e_2'=Te_1+e_2$,
$e_3'=e_3$.

{\em Case 2bb:} $\gamma_{2,3}\neq \gamma_{3,2}$. By setting $e_1'=e_1$,
$e_2'=e_2$, $e_3'= \tfrac{1}{\gamma_{2,3}-\gamma_{3,2}} e_3$ we see that we may
assume that $\gamma_{2,3}-\gamma_{3,2} = 1$. If $\gamma_{3,3}\neq -1$ then we
set $e_1'=e_1+e_3$, $e_2'=-\gamma_{2,3}e_1+e_2$, $e_3'=-\gamma_{3,3}e_1+e_3$,
and get $\theta(e_2',e_1')=\theta(e_3',e_1')=\theta(e_1',e_2')=\theta(e_1',e_3')
=0$, and $\theta(e_1',e_1')=\gamma_{3,3}+1\neq 0$. So here we are back in
Case 2a. If $\gamma_{3,3}=-1$ then set $e_1'=-e_1$, $e_2'=-e_2-\gamma_{2,3}e_3$,
$e_3'=e_3$, from which it is seen that we may assume that $\gamma_{2,3}=0$,
$\gamma_{3,2}=1$ as well. If $\gamma_{2,2}\neq 0$ then set $e_1'=e_2$,
$e_2'=e_2-\gamma_{2,2}e_1$, $e_3'=e_1-e_3$. The matrix of $\theta$ with respect
to this basis is $\gamma_{2,2}$ times
$$\left( \begin{smallmatrix} 1 & 1 & 0 \\
  0 & \gamma_{2,2} & -1 \\ 0 & -1 & 0 \end{smallmatrix} \right),$$
so that we are back in Case 2ba. If $\gamma_{2,2}=0$ then we obtain
$$\theta^5=\left( \begin{smallmatrix} 1 & 1 & 0 \\
  0 & 0 & 0 \\ 0 & 1 & -1 \end{smallmatrix} \right),$$
yielding $A_{4,11}$. Gr\"obner basis computations show that $\theta^5$
is not conjugate to the cocycles seen before.

\subsection{1-dimensional central extensions of $A_{3,2}$}

We have that $H^2(A_{3,2},F)$ consists of 
$\theta=\alpha_1 (\Delta_{a,b}+\Delta_{b,a})+\alpha_2\Delta_{a,c}+\alpha_3
\Delta_{c,a}+\alpha_4\Delta_{c,c}$. Furthermore
the automorphism group consists of 
$$\phi = \begin{pmatrix}
a_{11} & 0 & 0 \\ a_{21} & a_{11}^2 & a_{23}\\ a_{31} & 0 & a_{33} 
\end{pmatrix}.$$
Writing $\phi\theta =\alpha_1' (\Delta_{a,b}+\Delta_{b,a})+\alpha_2'\Delta_{a,c}+
\alpha_3'\Delta_{c,a}+\alpha_4'\Delta_{c,c}$ we have
\begin{align*}
\alpha_1' &= a_{11}^3 \alpha_1\\
\alpha_2' &= a_{11}a_{23}\alpha_1 +a_{11}a_{33}\alpha_2+a_{31}a_{33}\alpha_4\\
\alpha_3' &= a_{11}a_{23}\alpha_1 +a_{11}a_{33}\alpha_3+a_{31}a_{33}\alpha_4\\
\alpha_4' &= a_{33}^2\alpha_4.
\end{align*}

We need $\alpha_1\neq 0$ and ($\alpha_4\neq 0$ or $\alpha_2\neq \alpha_3$)
in order to have $\theta^\perp \cap C(A_{3,2})=0$. So after
dividing we may asume $\alpha_1=1$. Choose $a_{31}=0$, $a_{11}=1$ and
$a_{23}=-\alpha_3$. Then $\alpha_1'=1$, $\alpha_3'= 0$. So we may assume
$\alpha_3=0$.

First suppose that $\alpha_4\neq 0$.
Setting $a_{11}=a_{33}=\alpha_4$, and the other $a_{ij}$ equal to 0,
we obtain $\alpha_1'=\alpha_4'=
\alpha_4^3$, $\alpha_3'=0$. After dividing by $\alpha_4^3$ we see that we
may assume that $\alpha_4=1$ as well. If $\alpha_2\neq 0$ then we set $a_{22}
=\alpha_2^2$, $a_{33}=\alpha_2^3$, and the other $a_{ij}$ equal to 0,
leading to $\alpha_1'=\alpha_2'=\alpha_4'=\alpha_2^3$. Again, after dividing,
we conclude that we may assume that $\alpha_1=\alpha_2=\alpha_4=1$. So we get
two cocycles, $\Delta_{a,b}+\Delta_{b,a}+\Delta_{a,c}+\Delta_{c,c}$,
$\Delta_{a,b}+\Delta_{b,a}+\Delta_{c,c}$, yielding the algebras 
$A_{4,12}$, $A_{4,13}$. 
These are not conjugate, as one is symmetric and
the other is not.

Second, suppose that $\alpha_4=0$. Then $\alpha_2\neq \alpha_3$ implies that
$\alpha_2\neq 0$. Set $a_{23}=0$, $a_{11}=1$, $a_{33}=\tfrac{1}{\alpha_2}$,
showing that $\theta$ is conjugate to $\Delta_{a,b}+\Delta_{b,a}+\Delta_{a,c}$.
It is not conjugate to the previous ones, as cocycles with $\alpha_4\neq 0$
are not conjugate to cocycles with $\alpha_4=0$. This leads to the algebra
$A_{4,14}$. 

\subsection{1-dimensional central extensions of $A_{3,6}$}

Here $H^2(A_{3,6},F)$ is spanned by 
$\Delta_{b,b}+\Delta_{a,c}+\Delta_{c,a}$. So in this case we get only one algebra,
$A_{4,15}$.

\subsection{2-dimensional central extensions of $A_{2,1}$}
\label{sec:noncom2d}

Let $H= H^2(A_{2,1},F)$ which consists of all linear maps $A_{2,1}\to F$.
It is straightforward to see that every 2-dimensional
subspace of $H$ is usable. Therefore
the 2-dimensional central extensions of $A_{2,1}$ are parametrized by the
2-dimensional subspaces of $H$.

Let $a,b$ be a fixed basis of $A_{2,1}$. Then $\Delta_{a,a}$, $\Delta_{a,b}$,
$\Delta_{b,a}$, $\Delta_{b,b}$ form a basis of $H$. The 2-dimensional subspaces
are identified in the usual way with the points of a Grassmannian in
$\P(H\wedge H)$ (cf., \cite{shafI}, \S I.4.1). In $H\wedge H$ we use the
basis
$$\Delta_{a,a}\wedge \Delta_{a,b}, \Delta_{a,a}\wedge \Delta_{b,a},
\Delta_{a,a}\wedge \Delta_{b,b}, \Delta_{a,b}\wedge \Delta_{b,a},
\Delta_{a,b}\wedge \Delta_{b,b}, \Delta_{b,a}\wedge \Delta_{b,b}$$
(in that order). We write the homogeneous coordinates of a point
in $\P(V\wedge V)$, with respect to that basis, as $[\alpha_1,\ldots,
\alpha_6]$. By mapping the subspace with basis
$\theta_1,\theta_2\in H$ to the point $\theta_1 \wedge \theta_2\in
\P(H\wedge H)$, we obtain a bijection from the set of 2-dimensional
subspaces to the variety $\mathcal{X}$ of points
$[\alpha_1,\ldots,\alpha_6]\in \P(H\wedge H)$ with $\alpha_1\alpha_6-
\alpha_2\alpha_5+\alpha_3\alpha_4=0$. 

We have that $\Aut(A_{2,1})= \GL(A_{2,1})$.
Moreover, $\Aut(A_{2,1})$ acts on $H$ (see Section \ref{sec:method}), and
hence on $\mathcal{X}$. Moreover, by Proposition \ref{prop:2}, the
isomorphism classes of 2-dimensional central extensions of $A_{2,1}$ correspond
bijectively to the orbits of $\Aut(A_{2,1})$ on $\mathcal{X}$. 

As in Section \ref{sec4} we write an element of $\Aut(A_{2,1})$ as
$\phi=\left(\begin{smallmatrix} u & x \\ v & y \end{smallmatrix}\right)$ with
$uy-vx \neq 0$. Let $\alpha=[\alpha_1,\ldots,\alpha_6]\in \P(H\wedge H)$, then
$\phi(\alpha)=
[\beta_1,\ldots,\beta_6]$ with
\begin{equation}\label{eqn:conj0}
\begin{aligned}
  \beta_1 &= u^2 \alpha_1 +uv \alpha_3 -uv \alpha_4+v^2\alpha_6\\
  \beta_2 &= u^2 \alpha_2 +uv \alpha_3 +uv\alpha_4+v^2\alpha_5\\
  \beta_3 &= ux \alpha_1 + ux \alpha_2 +(uy+vx)\alpha_3 +vy \alpha_5 +vy \alpha_6\\
  \beta_4 &= -ux \alpha_1 +ux \alpha_2 +(uy+vx)\alpha_4 +vy \alpha_5 -vy \alpha_6\\
  \beta_5 &= x^2 \alpha_2+ xy \alpha_3+xy \alpha_4+y^2\alpha_5\\
  \beta_6 &= x^2 \alpha_1+xy\alpha_3-xy\alpha_4+y^2\alpha_6.
\end{aligned}
\end{equation}

Let $M$ be a subspace of $H$ with basis $\theta_1$, $\theta_2$. If
$\theta_i(m,m)=0$ for all $m\in A_{2,1}$, $i=1,2$, then both $\theta_i$ are
equal to a scalar multiple of $\Delta_{a,b}-\Delta_{b,a}$, and therefore
cannot be linearly independent. It follows that we may assume that
$\theta_1(a,a)=1$, and after subtracting a scalar multiple of $\theta_1$
from $\theta_2$, that $\theta_2(a,a)=0$. Represent the elements of $H$ by
their matrices with respect to the basis $a,b$ of $A_{2,1}$. Let $X$ be the
set of 2-dimensional subspaces with basis $\theta_1 = \left( \begin{smallmatrix}
  1 & \alpha \\ \beta & \gamma \end{smallmatrix}\right)$,
$\theta_2 = \left(\begin{smallmatrix} 0 & \delta \\ \epsilon & \eta
\end{smallmatrix}\right)$. We have just seen that every 2-dimensional
subspace of $H$ has a $\Aut(A_{2,1})$-conjugate in $X$. Furthermore, the
basis $\theta_1,\theta_2$ as above, corresponds to the point in
$\P(H\wedge H)$ with coordinates
\begin{equation}\label{eqcoord}
[ \delta, \epsilon, \eta, \alpha\epsilon-\beta\delta, \alpha\eta-
  \gamma\delta, \beta\eta-\gamma\epsilon].
\end{equation}
By $\widehat{X}$ we denote the image of $X$ in $\mathcal{X}$.
Then $\widehat{X}$ is exactly the set of points
$\alpha=[\alpha_1,\ldots,\alpha_6]\in\mathcal{X}$
with $(\alpha_1,\alpha_2,\alpha_3)\neq (0,0,0)$. We use \eqref{eqcoord} to
translate a point of $\widehat{X}$ to an algebra. In this section we only
deal with points in $\widehat{X}$, whereas in Section \ref{sec:list} we only
have algebras.

\begin{lem}\label{lem:A21.1}
  Set $\widehat{X}_1 = \{\alpha\in \widehat{X} \mid \alpha_1\neq 0\}$. Then
  an $\alpha \in \widehat{X}$ is conjugate to an element of $\widehat{X}_1$,
  or to $[0,1,0,0,0,0]$. Moreover, the latter point is not conjugate to an
  element of $\widehat{X}_1$.
\end{lem}

\begin{proof}
Let $\alpha \in \widehat{X}$,
and suppose that no $\Aut(A_{2,1})$-conjugate of $\alpha$ has first coordinate
nonzero. In particular $\alpha_1=0$ and the first coordinate of $\phi(\alpha)$
is $uv(\alpha_3-\alpha_4)+v^2\alpha_6$. Therefore, $\alpha_3=\alpha_4$ and
$\alpha_6=0$. It follows that $\alpha_2\alpha_5 = \alpha_3^2$. Hence
$\alpha_2\neq 0$ as otherwise $\alpha\not\in \widehat{X}$. It follows that
$\alpha = [0,1,\xi,\xi,\xi^2,0]$. Choose $u=1$, $v=0$, $y=1$,
$x=-\xi$. Then $\phi(\alpha)=[0,1,0,0,0,0]$. The conclusion is that if
$\alpha\in \widehat{X}$ has the property that none of its conjugates has
first coordinate nonzero, then $\alpha$ is conjugate to $[0,1,0,0,0,0]$.
The last statement is obvious from \eqref{eqn:conj0}.
\end{proof}

In the remainder of this section we study the orbits with representatives in
$\widehat{X}_1$. Since we may divide the homogeneous coordinates of a point
in $\P(V\wedge V)$ by a nonzero scalar, we may assume that the first coordinate
of a point in $\widehat{X}_1$ has first coordinate equal to 1.

\begin{lem}\label{lem:A21.2}
  Set $\widehat{X}_{1,0} = \{\alpha\in \widehat{X}_1 \mid \alpha_3= 0\}$.
  Let $\alpha\in \widehat{X}_1$. Then $\alpha$ is conjugate to an element of
  $\widehat{X}_{1,0}$ unless the characteristic of $F$ is $2$ and $\alpha_1=
  \alpha_2$, $\alpha_3\neq 0$ and $\alpha_4=0$, in which case $\alpha$ is
  conjugate to $p_\delta = [1,1,1,0,\delta,\delta]$. A point $p_\delta$ is
  not conjugate to points of $\widehat{X}_{1,0}$, and $p_\delta$, $p_\epsilon$
  are conjugate if and only if there is a $T\in F$ with $T^2+T+\delta+\epsilon
  =0$.
\end{lem}

\begin{proof}  
Let $\alpha= [\alpha_1=1,\alpha_2,\ldots,\alpha_6]\in \widehat{X}_1$.
Since $\alpha\in \mathcal{X}$, $\alpha_6 = \alpha_2\alpha_5-
\alpha_3\alpha_4$. Suppose $\alpha_3\neq 0$. Write $\phi(\alpha) = [\beta_1,
  \ldots,\beta_6]$ as in \eqref{eqn:conj0}.
Then $\beta_1 = u^2 +uv \alpha_3 -uv \alpha_4+
v^2(\alpha_2\alpha_5-\alpha_3\alpha_4)$ and $\beta_3 = ux(1+\alpha_2) +
(uy+vx)\alpha_3 +vy(\alpha_5 + \alpha_2\alpha_5-\alpha_3\alpha_4)$. If
$\alpha_2\neq -1$ then choose $v=0$, $u=y=1$,
$x=-\tfrac{1}{1+\alpha_2}\alpha_3$, so that $\beta_1=1$, $\beta_3=0$. Suppose
$\alpha_2=-1$, $\alpha_4\neq 0$. 
Then $\beta_3 = (uy+vx-vy\alpha_4)\alpha_3$. If $\alpha_5\neq 0$ then set
$x=0$, $y=u=1$, $v=\tfrac{1}{\alpha_4}$. Then $\beta_3=0$ and
$\beta_1 =-\alpha_5\alpha_4^{-2}\neq 0$. 
If $\alpha_5=0$ then set $u=0$, $v=x=1$,
$y=\tfrac{1}{\alpha_4}$. Then $\beta_3=0$, $\beta_1= -\alpha_3\alpha_4$ which
is nonzero. If $\alpha_2=-1$ and $\alpha_4=0$ then we have to choose
$u,v,x,y$ such that $uy+vx=0$, $uy-vx\neq 0$, $\beta_1= u^2+uv\alpha_3-v^2
\alpha_5\neq 0$. If the characteristic is not 2, then this clearly can be
done. Indeed, set $v=1$, $x=\tfrac{1}{2}$, $u\neq 0$ a non-zero of
$u^2+\alpha_3 u -\alpha_5$ (note that since we assume $\alpha_3\neq 0$
such a $u$ always exists even if $F$ is the field of three elements),
$y=-\tfrac{1}{2u}$. If the characteristic is 2, then if $\alpha_3\neq 0$,
$\alpha_2=1$, $\alpha_4=0$, 
by choosing $u=1$, $v= 0$, $y=\tfrac{1}{\alpha_3}$ we obtain $\beta_1=\beta_2
=\beta_3=1$, $\beta_4=0$.

By \eqref{eqn:conj0}, the polynomial
equations equivalent to $\phi p_\delta = \lambda p_\epsilon$ amount to
\begin{align*}
& u^2 +uv+v^2\delta+uy+vx =0\\
& x^2+xy+y^2\delta+(uy+vx)\epsilon=0\\
& uy+vx\neq 0.
\end{align*}
We claim that the above equations have a solution over $F$
if and only if there is a $T\in F$ with $T^2+t+\delta+\epsilon=0$
Indeed, if we have such a $T$ then we set $v=0$, $x=T$ and $u=y=1$ and obtain
a solution.
Conversely, the (reduced) Gr\"obner basis of the ideal generated by the above
polynomials (where we replace the last inequality by the polynomial
$D(uy+vx)+1$) contains the polynomials
\begin{align*}
  & (u+y)^2+v(u+y)+v^2(\delta+
  \epsilon)\\
  & v^2\epsilon^2 + vy\epsilon + x^2 + xy + y^2(\delta + \epsilon).
\end{align*}

So if a solution exists then those polynomials have to vanish as well.
If the solution has $v\neq 0$ then we divide the first polynomial by
$v^2$ and get $T=\tfrac{u+y}{v}$. If $v=0$ then $y\neq 0$ and from the
second polynomial we find $T= \tfrac{x}{y}$.
\end{proof}

Let $\alpha\in \widehat{X}_{1,0}$; we may assume that $\alpha_1=1$,
and hence $\alpha_6=\alpha_2\alpha_5$. Write $\beta = \phi(\alpha)$ as above.
If $\alpha_2=-1$ then $\beta_2 = -\beta_1$, so that, after dividing by $\beta_1$
we also have $\beta_2=-1$. Therefore an $\alpha\in \widehat{X}_{1,0}$ with
$\alpha_1=1$ and $\alpha_2=-1$ is not conjugate to 
$\alpha'\in \widehat{X}_{1,0}$ with $\alpha_1'=1$ 
$\alpha_2'\neq -1$.

\begin{lem}\label{lem:A21.3A}
  Let the characteristic be different from 2.
  Let $\alpha \in \widehat{X}_{1,0}$ with $\alpha_1=1$, $\alpha_2=-1$.
  Then $\alpha$ is conjugate to $q_\delta = [1,-1,0,0,\delta,-\delta]$. Moreover,
  $q_\delta$, $q_\epsilon$ are conjugate if and only if there is a $\nu\in F^*$
  with $\epsilon = \nu^2\delta$.
\end{lem}  

\begin{proof}
Again write $\beta = \phi(\alpha)$, as in \eqref{eqn:conj0}. 
Then $\beta_3=0$. Furthermore, $\beta_4=
-2ux +(uy+vx)\alpha_4 +2vy \alpha_5$. By taking $v=0$, $u=y=1$,
$x=\tfrac{1}{2}\alpha_4$ we obtain $\beta_4=0$.
So $\alpha$ is conjugate to $q_\delta$. Write $\beta=\phi(q_\delta)$. Then
$\beta_1 = u^2-\delta v^2$, $\beta_2=-\beta_1$, $\beta_3=0$,
$\beta_4=-2ux+2\delta vy$, $\beta_5 = -x^2+\delta y^2$. Hence $p_\delta$
is conjugate to $p_\epsilon$ if and only if
there are $u,v,x,y$ with $uy-vx\neq 0$, $u^2-\delta v^2\neq 0$,
$ux-\delta vy=0$, $-x^2+\delta y^2=( u^2-\delta v^2)\epsilon$.
A Gr\"obner basis of the ideal generated by these polynomials contains
$u^2\epsilon - y^2\delta$,  $v^2\delta\epsilon - x^2$. It follows that if
$p_\delta$, $p_\epsilon$ are conjugate, then there is a nonzero $\nu\in F$ with
$\epsilon = \nu^2 \delta$. The converse is straightforward, by setting
$v=x=0$, $y=1$, $u=\tfrac{1}{\nu}$.
\end{proof}

\begin{lem}\label{lem:A21.3B}
  Let the base field have characteristic 2.
  Let $\alpha \in \widehat{X}_{1,0}$ with $\alpha_1=\alpha_2=1$.
  If $\alpha_4=0$, then
  $\alpha=q_\delta = [1,1,0,0,\delta,\delta]$. Moreover,
  $q_\delta$, $q_\epsilon$ are conjugate if and only if there are $u,v,x,y\in F$
  with $uy+vx\neq 0$, $u^2+v^2\delta\neq 0$ and
  $\epsilon = \tfrac{x^2+y^2\delta}{u^2+v^2\delta}$.
  If $\alpha_4\neq 0$
  then $\alpha$ is not conjugate to a $q_\delta$ but to $r_\delta = [1,1,0,1,
  \delta,\delta]$. Finally, $r_\delta$, $r_\epsilon$ are conjugate if and only
  if there is a $T\in F$ with $T^2+T+\delta+\epsilon=0$. 
\end{lem}

\begin{proof}
  The first statement is obvious. The conjugacy condition follows directly from
  the polynomials already written in the proof of the previous lemma.
  
Suppose that $\alpha_4\neq 0$, then in \eqref{eqn:conj0} we have $\beta_4=
(uy+vx)\alpha_4 \neq 0$. Therefore $\alpha$ is not conjugate to a $q_\delta$.
In \eqref{eqn:conj0}
we take $v=0$, $u=1$, $y=\tfrac{1}{\alpha_4}$ and obtain $\beta_4=1$. So in
this case $\alpha$ is conjugate to $r_\delta$. By
\eqref{eqn:conj0}, $r_\delta$ is conjugate to $r_\epsilon$ if and only if there
are $u,v,x,y\in F$ with $uy+vx\neq 0$,
$uy+vx=u^2+uv+\delta v^2$, $x^2+xy+\delta y^2=(u^2+uv+\delta v^2)\epsilon$.
By a Gr\"obner basis computation it is seen that this implies that
$v^2(\delta+\epsilon)+v(u+y)+(u+y)^2=0$, $v^2\epsilon^2+vy\delta+x^2+xy
+y^2\delta+y^2\epsilon=0$. As in the proof of Lemma \ref{lem:A21.2}, 
this implies that there is a $T\in F$ with $T^2+T+
\delta+\epsilon=0$. Conversely, if such a $T$ exists,
then $p_\delta$, $p_\epsilon$ are seen to be conjugate by setting $v=0$, $u=y=1$,
$x=T$. 
\end{proof}

\begin{lem}\label{lem:A21.4}
  Set $\widehat{Y}= \{\alpha\in \widehat{X}_{1,0} \mid \alpha_1=1, \alpha_2
  \neq -1\}$. Let $\alpha,\alpha'\in \widehat{Y}$ be conjugate. Then there
  is a $\nu\in F^*$ with $\alpha_5' = \nu^2\alpha_5$. Conversely $\alpha\in
  \widehat{Y}$ is
  conjugate to 
  $[1,\alpha_2,0,\nu\alpha_4,\nu^2\alpha_5,\nu^2\alpha_6]\in
  \widehat{Y}$.
\end{lem}  

\begin{proof}
Write $\beta = \phi(\alpha)$ as in \eqref{eqn:conj0}. 
Then $\beta_3 = (1+\alpha_2)(ux+\alpha_5
vy)$. Furthermore, $\beta_1 =  u^2 -uv \alpha_4+v^2\alpha_2\alpha_5$,
$\beta_5=x^2 \alpha_2+ xy \alpha_4+y^2\alpha_5$. Suppose that $\alpha$ is
conjugate to $\alpha'=[1,\alpha_2',0,\alpha_4',\alpha_5',\alpha_2'\alpha_5']$.
Then there are $u,v,x,y\in F$ with $uy-xv\neq 0$, $u^2 -uv \alpha_4+
v^2\alpha_2\alpha_5\neq 0$, $ux+\alpha_5 vy=0$, $x^2 \alpha_2+ xy \alpha_4+
y^2\alpha_5 = (u^2 -uv \alpha_4+v^2\alpha_2\alpha_5)\alpha_5'$. (Of course,
there are further conditions coming from $\beta_2$, $\beta_4$, but we ignore
those.) By a Gr\"obner basis computation it is seen that
these equations imply $u^2\alpha_5' - y^2\alpha_5=0$, $v^2\alpha_5\alpha_5'
-x^2=0$. But that implies that there is a nonzero $\nu\in F$ with
$\alpha_5' = \nu^2 \alpha_5$. For the converse set 
$v=x=0$, $u=1$, $y=\nu$.
\end{proof}

Let $\widehat{Z}_0$, $\widehat{Z}_1$ be the sets of points of $\widehat{Y}$
with, respectively, fifth coordinate zero, and fifth coordinate nonzero.
Then in particular it follows that
points from $\widehat{Z}_0$ are not conjugate to points from $\widehat{Z}_1$.

\begin{lem}\label{lem:A21.5}
  Let $\alpha\in \widehat{Z}_0$. Then $\alpha$ is conjugate to 
  $[1,0,0,-1,0,0]$ or to $[1,\delta,0,0,0,0]$, $\delta\in F$, $\delta\neq 0$.
  These points are pairwise not conjugate.
\end{lem}

\begin{proof}
Suppose that $\alpha_4=0$
as well, and write $\beta=\phi(\alpha)$ as in \eqref{eqn:conj0}.
Then $\beta_1 = u^2$,
$\beta_2 = u^2\alpha_2$, $\beta_3 = (1+\alpha_2)ux$, $\beta_4=ux(-1+\alpha_2)$,
$\beta_5=x^2 \alpha_2$. So if $[\beta_1,\ldots,\beta_6]$ lies in $\widehat{Y}$,
then $u\neq 0$ and $x=0$ so that also $\beta_4=0$. It follows that points from
$\widehat{Z}_0$ with fourth coordinate zero are not conjugate to points
of $\widehat{Z}_0$ with fourth coordinate nonzero. Moreover, it follows that
$[1,\delta,0,0,0,0]$ and $[1,\epsilon,0,0,0,0]$ are
conjugate if and only if $\delta=\epsilon$. Second,
if $\alpha_4\neq 0$ then set $x=0$, $u=-\alpha_4$, $v=\alpha_2$, $y=1+\alpha_2$
and see that $\alpha$ is conjugate to $[1,0,0,-1,0,0]$.
\end{proof}

Let $\alpha\in \widehat{Z}_1$, and set $\beta=\phi(\alpha)$ as before. Then
\begin{equation}\label{eqn:conj}
\begin{aligned}
  \beta_1 &= u^2  -uv \alpha_4+v^2\alpha_2\alpha_5\\
  \beta_2 &= u^2 \alpha_2 +uv\alpha_4+v^2\alpha_5\\
  \beta_3 &= (1+\alpha_2)(ux +vy \alpha_5)\\
  \beta_4 &= ux(-1 +\alpha_2) +(uy+vx)\alpha_4 +vy(1-\alpha_2)\alpha_5\\
  \beta_5 &= x^2 \alpha_2+ xy \alpha_4+y^2\alpha_5
\end{aligned}
\end{equation}

\begin{lem}\label{lem:A21.6a}
  Let $\alpha,\alpha'\in \widehat{Z}_1$ be such that 
  $\alpha_5=\alpha_5'$.
Then $\alpha,\alpha'$ are conjugate if and only if there are $u,v\in F$,
$\epsilon = \pm 1$ with $u^2+v^2\alpha_5 \neq 0$, $u^2  -uv \alpha_4+v^2\alpha_2
\alpha_5 \neq 0$,
and after setting $y=\epsilon u$, $x=-\epsilon v\alpha_5$ we have
$\alpha_i' = \tfrac{\beta_i}{\beta_1}$, where the $\beta_i$ are as in
\eqref{eqn:conj}. 
\end{lem}  

\begin{proof}
  Suppose that $\alpha,\alpha'$ are conjugate. Write $\beta=\phi(\alpha)$ as
  in \eqref{eqn:conj} and suppose $\beta=\alpha'$. 
As seen in the proof of Lemma \ref{lem:A21.4}, 
$uy-vx\neq 0$, $u^2  -uv \alpha_4+v^2\alpha_2
\alpha_5 \neq 0$ and 
$u^2\alpha_5 - y^2\alpha_5=0$,  $v^2\alpha_5^2 - x^2=0$. Hence $y=\epsilon u$,
$x=\nu v\alpha_5$, with $\epsilon,\nu = \pm 1$. Also we must have
$ux +vy \alpha_5=0$. If $uv\neq 0$ then this yields $\nu=-\epsilon$. On the
other hand, if $u=0$ then $y=0$ and we can choose $\epsilon=-\nu$. Similarly,
if $v=0$ then $x=0$ and we can choose $\nu = -\epsilon$. Finally, $uy-vx\neq 0$
translates to $u^2+v^2\alpha_5\neq 0$. The other direction
is trivial.
\end{proof}

\begin{lem}\label{lem:A21.6}
  Let $\alpha\in \widehat{Z}_1$. If $\alpha_4\neq 0$ 
  then $\alpha$ is conjugate to an $\alpha'
  \in \widehat{Z}_1$ with $\alpha_2'\neq 1$, unless
  $F$ has three elements and $\alpha_5=-1$, in which case $\alpha$ is
  conjugate to $[1,1,0,-1,-1,-1]$. On the other hand, if $\alpha_2=1$ and
  $\alpha_4=0$ then the characteristic is not $2$, $\alpha$ is not conjugate
  to an $\alpha'\in \widehat{Z}_1$ with $\alpha_2'\neq 1$, but to $s_\delta =
  [1,1,0,0,\delta,\delta]$. We have that $s_\delta$ is conjugate to $s_\epsilon$
  if and only if there is a $\nu\in F^*$ with $\epsilon = \nu^2 \delta$.
  \end{lem}

\begin{proof}
Suppose that $\alpha_2=1$, $\alpha_4\neq 0$.
(By hypothesis $\alpha_2\neq -1$, so in particular the
characteristic is not 2.) By the previous lemma $\alpha$ is conjugate to an
$\alpha'\in
\widehat{Z}_1$ with $\alpha_2'\neq 1$ if and only if there are $u,v\in F$
with $-uv\alpha_4 \neq uv\alpha_4$ (this follows from $\beta_1\neq \beta_2$ 
in \eqref{eqn:conj}),  $0\neq u^2+v^2\alpha_5$,
$0\neq u^2  -uv \alpha_4+v^2\alpha_5$.
Choose $v=1$ and $u\neq 0$ such that $u^2\neq
-\alpha_5$, $u^2-\alpha_4 u +\alpha_5\neq 0$. If $F$ has more than five
elements then such $u$ clearly exist. If $F$ has five elements, then it is
not possible that two nonzero elements of $F$ are solutions of
$X^2=-\alpha_5$, and two other nonzero elements are solutions of
$X^2-\alpha_4X+\alpha_5=0$, because $\alpha_4\neq 0$; so also in that case
a $u$ as above exists. If $F$ has three elements,
then the second equation cannot have two distinct roots as $\alpha_4\neq 0$
and the first equation has roots only if $\alpha_5=-1$. So if $\alpha_5\neq -1$
then we can find a $u$ as above. If $\alpha_5=-1$ then as $|F|=3$, $u=\pm v$
and $u^2+v^2\alpha_5=0$. Hence $\alpha$ is not conjugate to a point with
second coordinate $\neq 1$. In this case, if $\alpha_4=1$ then by
\eqref{eqn:conj} with $v=x=0$, $u=1$, $y=-1$ we see that
$\alpha$ is conjugate to $[1,1,0,-1,-1,-1]$. We conclude that, if
$\alpha_4\neq 0$ then $\alpha$ is conjugate to $\alpha'\in \widehat{Z}_1$
with $\alpha_2'\neq 1$, unless $|F|=3$, in which case there is an
extra point.

If $\alpha_2=1$ and $\alpha_4=0$, then $\alpha = s_\delta$.
The conjugacy condition is seen in Lemma \ref{lem:A21.4}.
\end{proof}

Let $\widehat{W}$ denote the set of $\alpha\in \widehat{Z}_1$ with $\alpha_1=1$,
$\alpha_2\neq 1$. For $\alpha\in \widehat{W}$ define
$$\sigma(\alpha) = \frac{\alpha_4}{1-\alpha_2}.$$
Let $\alpha \in \widehat{W}$. Then also 
$\alpha_3=0$, $\alpha_5\neq 0$.
Let $\alpha'\in \widehat{W}$ be such that 
$\alpha_5'=\alpha_5$. By Lemma \ref{lem:A21.6a}, $\alpha$, $\alpha'$ are
conjugate if and only if there are $u,v\in F$ with $u^2+\alpha_5 v^2\neq 0$,
$u^2  -uv \alpha_4+v^2\alpha_2\alpha_5\neq 0$ and
\begin{equation}\label{eqn:conj2}
\begin{aligned}
  \alpha_2' &= \frac{u^2 \alpha_2 +uv\alpha_4+v^2\alpha_5}
        {u^2  -uv \alpha_4+v^2\alpha_2\alpha_5}\\
  \alpha_4' &= \epsilon \frac{2uv\alpha_5(1-\alpha_2)+(u^2-v^2\alpha_5)\alpha_4}
                {u^2  -uv \alpha_4+v^2\alpha_2\alpha_5}.
\end{aligned}
\end{equation}
Brute force verification shows that \eqref{eqn:conj2} implies 
\begin{equation}\label{eq:sigma}
  \sigma(\alpha')^2+\alpha_5 = \Psi^2 (\sigma(\alpha)^2+\alpha_5)
  \text{ with } \Psi=
    \frac{(\alpha_2-1)(u^2+v^2\alpha_5)}
         {(\alpha_2-1)(v^2\alpha_5-u^2)-2uv\alpha_4}.
\end{equation}         

\begin{lem}\label{lem:A21.7}
Suppose that the characteristic of $F$ is not 2.
Let $\alpha,\alpha'\in\widehat{W}$.
Suppose that $\sigma(\alpha)^2=-\alpha_5$. Then $\alpha$ is conjugate to
$[1,0,0,1,-1,0]$. Furthermore, $\alpha$, $\alpha'$ are conjugate if and
only if $\sigma(\alpha')^2 = -\alpha_5'$. 
\end{lem}

\begin{proof}
Choosing $v=x=0$, $y=1$ and
$u=\sigma(\alpha)$ we obtain by \eqref{eqn:conj} that $\beta_1=u^2$,
$\beta_2=u^2\alpha_2$, $\beta_3=0$, $\beta_4 = u\alpha_4$, $\beta_5 =
\alpha_5$. Dividing by $u^2$ (note that $\sigma(\alpha)\neq 0$ by hypothesis),
we see that $\alpha$ is conjugate to $[1,\alpha_2,0,1-\alpha_2,-1,0]$.
Now we use the formulas
\eqref{eqn:conj2} (where instead of $\alpha_4$ we put $1-\alpha_2$, instead
of $\alpha_5$ we put $-1$).
Setting $u=-1$, $v=\alpha_2$, $\epsilon=1$, we have
$u^2+\alpha_5 v^2 = 1-\alpha_2^2\neq 0$,
and $u^2  -uv \alpha_4+v^2\alpha_2\alpha_5 = (1-\alpha_2)(1+\alpha_2)^2\neq 0$,
and $\alpha_2' = 0$, $\alpha_4'=1$, showing that $\alpha$ is conjugate to
$[1,0,0,1,-1,0]$.

Suppose $\alpha,\alpha'$ are conjugate. By Lemma \ref{lem:A21.4}, there is a
$\nu\in F^*$ such that $\alpha_5 = \nu^2\alpha_5'$, and moreover,
$\alpha'$ is conjugate to $[1,\alpha_2',0,\nu\alpha_4',\alpha_5,\alpha_6']$.
By \eqref{eq:sigma} this implies that $\nu^2\sigma(\alpha')^2 +\alpha_5 =
\Psi^2(\sigma(\alpha)^2 +\alpha_5)$. The denominator of $\Psi$ is
$\tfrac{(u(\alpha_2-1)+v\alpha_4)^2}{\alpha_2-1}$. Hence it is zero if and
only if $u=v\sigma(\alpha)$. But then $u^2+\alpha_5v^2=0$.
The conclusion is that 
necessarily $\nu^2\sigma(\alpha')^2=-\alpha_5$, or equivalently,
$\sigma(\alpha')^2 =-\alpha_5'$. The converse is obvious, as
$\sigma(\alpha')^2 =
-\alpha_5'$ implies that $\alpha'$ is conjugate to $[1,0,0,1,-1,0]$ as well.
\end{proof}

\begin{lem}\label{lem:A21.8A}
Suppose that the characteristic of $F$ is not 2.
Let $\alpha,\alpha'\in\widehat{W}$ be such that $\alpha_5=\alpha_5'$.
Suppose that $\sigma(\alpha)=\sigma(\alpha')$ and $\sigma(\alpha)^2\neq
-\alpha_5$. Then $\alpha,\alpha'$ are
conjugate if and only if $\alpha_2'=\alpha_2$, $\alpha_4'=\alpha_4$ or
$\alpha_2'=\alpha_2^{-1}$, $\alpha_4'=-\alpha_4\alpha_2^{-1}$.
\end{lem}

\begin{proof}
If $\alpha,\alpha'$ are conjugate, then there are $u,v\in F$ with
$u^2+\alpha_5 v^2\neq 0$, $u^2  -uv \alpha_4+v^2\alpha_2\alpha_5\neq 0$ and
\eqref{eqn:conj2}. Write $\sigma=\sigma(\alpha)$, $\sigma'=\sigma(\alpha')$.  
If $\epsilon=1$ then $\sigma'=\sigma$ amounts to
$uv(\sigma^2+\alpha_5)=0$ so that $uv=0$. If $u=0$ then $\alpha_2'=
\alpha_2^{-1}$, $\alpha_4'=-\alpha_4\alpha_2^{-1}$. If $v=0$ then $\alpha_2'=
\alpha_2$, $\alpha_4'=\alpha_4$. If $\epsilon=-1$ then $\sigma=\sigma'$
is equivalent to $(\alpha_4u+\alpha_5(1-\alpha_2)v)(\alpha_4 v -(1-\alpha_2)u)
=0$. If the first factor vanishes then $v=-\alpha_5^{-1}\sigma u$, and
$\alpha_2'=\alpha_2$, $\alpha_4'=\alpha_4$. If the second factor is zero then
$\alpha_2'=\alpha_2^{-1}$, $\alpha_4'=-\alpha_4\alpha_2^{-1}$.
\end{proof}  

\begin{lem}\label{lem:A21.8}
Suppose that the characteristic of $F$ is not 2.
Let $\alpha,\alpha'\in\widehat{W}$ and write $\sigma=\sigma(\alpha)$,
$\sigma'=\sigma(\alpha')$.
Suppose that $\sigma^2\neq -\alpha_5$, $(\sigma')^2\neq -\alpha_5'$.
If $\alpha,\alpha'$ are conjugate, then there are $\varphi,\nu\in F^*$ with
$(\sigma')^2+\alpha_5' = \varphi^2(\sigma^2+\alpha_5)$,
$\alpha_5 = \nu^2\alpha_5'$. Conversely, suppose that these conditions are
satisfied, set $\psi = \nu\varphi$,
let $\omega_1= \pm 1$ be such that $\psi-\omega_1
\neq \tfrac{2\psi}{1-\alpha_2}$, and set
$$\gamma_2 = \frac{\psi(1+\alpha_2)-\omega_1(1-\alpha_2)}{
  \psi(1+\alpha_2) + \omega_1(1-\alpha_2)}\text{ and }
\gamma_4= \nu\sigma'(1-\gamma_2).$$
Then $\alpha,\alpha'$ are conjugate if and only if $\alpha_2'=\gamma_2$,
$\nu\alpha_4' = \gamma_4$, or $\alpha_2'=\gamma_2^{-1}$, $\nu\alpha_4'=-\gamma_4
\gamma_2^{-1}$. 
\end{lem}

\begin{proof}
Suppose that $\alpha$, $\alpha'$ are conjugate. By Lemma \ref{lem:A21.4}
there is a $\nu \in F^*$ with $\alpha_5 = \nu^2 \alpha_5'$, and $\alpha'$
is conjugate to $[1,\alpha_2',0,\nu\alpha_4',\alpha_5,\alpha_6']$. So also
$\alpha$ is conjugate to this point. By \eqref{eq:sigma}
there is a $\psi\in F^*$
with $\nu^2 (\sigma')^2 +\alpha_5 = \psi^2(\sigma^2+\alpha_5)$. Dividing by
$\nu^2$ we see that $(\sigma')^2+\alpha_5' = \varphi^2(\sigma^2+\alpha_5)$
with $\varphi = \psi/\nu$.

For the converse, define $\omega_1$ as in the statement of the lemma. Set
$\tau = \nu\sigma'$. Then $\tau^2+\alpha_5 = \psi^2(\sigma^2+\alpha_5)$.

Suppose that $\psi=\omega_1$. Then $\tau = \pm \sigma$. If both are zero, then
$\alpha'$ is conjugate to $[1,\alpha_2',0,0,\alpha_5,\alpha_6']$. Furthermore,
$\alpha=[1,\alpha_2,0,0,\alpha_5,\alpha_6]$ and $\gamma_2=\alpha_2$,
$\gamma_4=0$. So Lemma \ref{lem:A21.8A} finishes the proof in this case.

If $\psi\neq \omega_1$ or both $\sigma,\tau$ are nonzero, then there is an
$\omega_2=\pm 1$ such that $\omega_1\psi(\omega_2\tau\sigma-
\alpha_5) \neq -\tau^2-\alpha_5$. 
Set $u=\omega_2\tau + \omega_1\psi
\sigma$, $v=\omega_1\psi-1$. After some manipulation it is seen that
\begin{align*}
  u^2-2uv\sigma -\alpha_5v^2 &= 2\omega_1\psi^{-1} (\tau^2+\omega_1
  \omega_2 \tau\psi\sigma-\omega_1\alpha_5\psi+\alpha_5)\\
  u^2+\alpha_5 v^2 &= 2(\tau^2+\omega_1\omega_2
  \tau\psi\sigma-\omega_1\alpha_5\psi+\alpha_5)\\
  u^2-\alpha_4uv+\alpha_2\alpha_5v^2 &= \tfrac{1}{2}\psi^{-1} (u^2+
  \alpha_5 v^2)(2\psi -(\psi-\omega_1)(1-\alpha_2))\\
  \alpha_2 u^2 +\alpha_4 uv + \alpha_5 v^2 &= \tfrac{1}{2}\psi^{-1} (u^2+
  \alpha_5 v^2)(2\alpha_2\psi +(\psi-\omega_1)(1-\alpha_2)).
\end{align*}
In particular we see that $u^2+\alpha_5 v^2\neq 0$ by the choice of $\omega_2$.
So these $u,v$ define an element $\phi$ of $\Aut(A_{2,1})$. By the choice
of $\omega_1$ we see that $u^2-\alpha_4 uv +\alpha_2\alpha_5 v^2\neq 0$. Write
$\gamma = \phi(\alpha)$. Then $\gamma_2,\gamma_4$ are given by the right
hand sides of \eqref{eqn:conj2}, and $\gamma_1=1$, $\gamma_3=0$, $\gamma_5=
\alpha_5$. Secondly, $\gamma_2$ is given as in the statement of the lemma.
Thirdly, the factor $\Psi$ in \eqref{eq:sigma} is equal
to $-\omega_1\psi$. Hence $\sigma(\gamma)^2+\alpha_5 = \psi^2(
\sigma^2+\alpha_5) = \tau^2+\alpha_5$, so that $\sigma(\gamma) = \pm \tau$.
By choosing $\epsilon$ in \eqref{eqn:conj2}
we can force $\sigma(\gamma)=\tau$. This ensures
that $\gamma_4 = \nu\sigma'(1-\gamma_2)$. Furthermore, $\gamma_2\neq 1$ as
otherwise $\omega_1=0$, so we have $\gamma\in\widehat{W}$. By Lemma
\ref{lem:A21.4}, $\alpha'$ is conjugate to $\delta = [1,\alpha_2',0,
\nu\alpha_4',\alpha_5,\delta_6]$. 
We have $\sigma(\delta) = \nu\sigma'$.
Since $\alpha$, $\alpha'$ are conjugate if and only if $\gamma$, $\delta$ are,
Lemma \ref{lem:A21.8A} finishes the proof.
\end{proof}

\begin{rem}
If both choices of $\omega_1=\pm 1$ are possible then it does not matter which
one is chosen. Indeed, if $\omega_1=1$ yields $\gamma_2$, $\gamma_4$, then
$\omega_1=-1$ yields $\gamma_2^{-1}$, $-\gamma_4\gamma_2^{-1}$. 
\end{rem}

\begin{lem}\label{lem:A21.9}
Suppose that the characteristic of $F$ is 2.
Let $\alpha,\alpha'\in\widehat{W}$ and write $\sigma=\sigma(\alpha)$,
$\sigma'=\sigma(\alpha')$. Define
$$H_{\alpha} = \{ \frac{\sigma uv+\alpha_5 v^2}{u^2+\alpha_5 v^2} \mid
u,v\in F \text{ and } u^2+\alpha_5 v^2\neq 0\}.$$
Then $H_{\alpha}$ is a subgroup of the additive group of $F$.
Moreover, $\alpha$, $\alpha'$ are conjugate if and only if there is
a $\nu\in F^*$ with $\alpha_5 = \nu^2 \alpha_5'$, $\sigma = \nu \sigma'$
and a $h\in H_\alpha$
such that $\tfrac{1}{1+\alpha_2'} = \tfrac{1}{1+\alpha_2} +h$. 
\end{lem}

\begin{proof}
By direct computation it is verified that $H_\alpha$ is a subgroup of $F$.   
Suppose that $\alpha$, $\alpha'$ are conjugate. By Lemma \ref{lem:A21.4}
there is a $\nu\in F^*$ with $\alpha_5 = \nu^2 \alpha_5'$ and  
$\alpha'$ is conjugate to $\delta = [1,\alpha_2',0,
\nu\alpha_4',\alpha_5,\delta_6]$. Let $\alpha$ be conjugate to $\gamma$ where
$\gamma_2$, $\gamma_4$ are the right hand sides of \eqref{eqn:conj2}, and
$\gamma_1=1$, $\gamma_3=0$, $\gamma_5=\alpha_5$. Because $\alpha$, $\alpha'$
are conjugate, $u$, $v$ can be chosen such that $\gamma=\delta$.
Since the characteristic is 2,
the $\Psi$ of \eqref{eq:sigma} is 1. Hence, by the same equation,
$\sigma(\alpha)=\sigma(\gamma)=\sigma(\delta) = \nu\sigma'$.

Now suppose that a $\nu$ satisfying the given conditions exists. Let $\delta$,
$\gamma$ be as above. If $\alpha$, $\alpha'$ are conjugate, $u,v$ can be chosen
such that $\gamma=\delta$. But
$$\frac{1}{1+\gamma_2} = \frac{u^2+uv\alpha_4+v^2\alpha_2\alpha_5}
{u^2+\alpha_5v^2}\cdot \frac{1}{1+\alpha_2} = \frac{1}{1+\alpha_2}+
\frac{\sigma uv+\alpha_5 v^2}{u^2+\alpha_5 v^2}.$$
Conversely, if $u,v\in F$ exist with $u^2+\alpha_5 v^2\neq 0$ and
$\tfrac{1}{1+\alpha_2'} = \frac{1}{1+\alpha_2}+
\frac{\sigma uv+\alpha_5 v^2}{u^2+\alpha_5 v^2}$, then $u^2+uv\alpha_4+
v^2\alpha_2\alpha_5\neq 0$ because $\tfrac{1}{1+\alpha_2'}
\neq 0$. So we can define $\phi$ with these $u,v$ and $\gamma = \phi(\alpha)$.
Then $\tfrac{1}{1+\alpha_2'} = \tfrac{1}{1+\gamma_2}$, so $\gamma_2=\alpha_2'$.
Furthermore, $\tfrac{\gamma_4}{1+\gamma_2}=\sigma(\gamma)=\sigma(\alpha)
=\nu \sigma' = \sigma(\delta) = \tfrac{\nu \alpha_4'}{1+\alpha_2'}$. Hence
$\gamma_4=\nu \alpha_4'$. It is seen that $\alpha$ is conjugate to $\delta$,
and therefore to $\alpha'$.
\end{proof}

\section{The classification over specific fields}\label{sec:fields}

Here we give complete and irredundant lists of the 4-dimensional nilpotent
associative algebras over finite fields and over $\R$. We comment on the
classification over algebraically closed fields (where we can get an
explicit list for example if $F=\C$)
and over $\Q$ (where we cannot obtain a fully explicit list).

\subsection{Finite fields of odd characteristic}

In this section $F$ is a finite field of odd characteristic.

\begin{lem}\label{lem:charodd1}
There are non-squares $\xi\in F$ such that $\xi-1$ is a square.
\end{lem}

\begin{proof}
  Let $p$ be the characteristic of $F$. Suppose that there are $i$ with
  $1 < i \leq p-1$ that are non-squares. Then let $\xi$ be the minimal such
  $i$. On the other hand, if there are no such $i$, then $-1=p-1$ is a
  square. Let $\eta\in F$ be a primitive element. If $\eta-1$ is a square
  then we can take $\xi=\eta$. If $\eta-1$ is not a square then $\eta-1 =
  \eta^{2k+1}$ and $\eta^{-1}-1 = -\eta^{2k}$ is a square, so $\xi = \eta^{-1}$
  does the job.
\end{proof}

\begin{lem}\label{lem:charodd2}
There are non-squares $\xi\in F$ such that $1-\xi$ is a square, unless
$|F|=3$ where there is no such $\xi$.  
\end{lem}

\begin{proof}
  Suppose that $-1$ is a square. Let $\xi\in F$ be a non-square such that
  $\xi-1$ is a square (previous lemma). Then $1-\xi = -(\xi-1)$ is a square
  as well.

  Suppose that $-1$ is not a square. If $2$ is a square then take $\xi=-1$.
  So suppose that $2$ is not a square. Suppose also that the characteristic
  is not 3. If $3$ is not a square then we can take $\xi = 3$ as $-2$ is a
  square. So suppose that $3$ is a square. But then $-3$ is not a square
  and $1-(-3)=4$ is a square.

  There remains the case where the characteristic is $3$ and $-1=2$ is not a
  square. Suppose that $F$ is such that there is no non-square $\xi$ such
  that $1-\xi$ is a square. Let $f : F\to F$ be the map with
  $f(\alpha) =1-\alpha$. Let $S$ be the set of squares in $F$, except $0,1$.
  Let $N$ be the set of non-squares in $F$, except $-1$. Suppose that $|F|>3$
  so that $S$, $N$ are non-empty.  Then $f$ maps
  $N$ to $N$, and hence also $S$ to $S$. If $\zeta\in S$ is such that
  $1+\zeta\in S$, then we $\xi = -\zeta$ is not a square with $1-\xi$
  a square. So under our assumption such $\zeta$ do not exist. Define
  $g : F\to F$, $g(\alpha) = 1+\alpha$. Then $g$ maps $S$ to $N$,
  and hence $N$ to $S$. Let $\zeta\in N$. Then $1-\zeta\in N$, $1+\zeta\in
  S$. But then $\zeta^2\in S$ and $f(\zeta^2)=(1-\zeta)(1+\zeta)\in N$,
  which is a contradiction.
\end{proof}

\begin{prop}
  Let $\eta$ be a primitive element of $F$. Let $\xi,\zeta\in F$ be non-squares
  such that $\xi-1$, $1-\zeta$ are squares. Fix $\sigma_\xi, \sigma_\zeta\in F$
  such that $\sigma_\xi^2 = \xi-1$, $\sigma_\zeta^2=1-\zeta$. Let $\mathcal{A}$
  be a maximal subset of $F$ with $0\in \mathcal{A}$,
  $\pm 1 \not\in \mathcal{A}$, if $\alpha\in \mathcal{A}$, $\alpha\neq 0$
  then $\alpha^{-1}\not\in \mathcal{A}$. Then
  \begin{align*}
    & A_{4,1}, A_{4,2}, A_{4,3}^\delta (\delta=1,\eta), A_{4,4}^\delta (\delta\in F),
    A_{4,5}, A_{4,6}, A_{4,7}, A_{4,8}^{1,1}, A_{4,9}^{1,\beta} (\beta\in F,
    \beta\neq \tfrac{1}{4}),  \\
    & A_{4,9}^{1,\tfrac{1}{4}},  A_{4,9}^{\eta,\tfrac{1}{4}}, A_{4,10}^0, A_{4,11},
    A_{4,12}, A_{4,13}, A_{4,14}, A_{4,15}, A_{4,16}, A_{4,18}^\delta
    (\delta = 0,1,\eta), A_{4,20}\\
    & A_{4,21}^\delta (\delta\in F, \delta\neq -1),
    A_{4,23}^\delta (\delta = 1,\eta), A_{4,24}, A_{4,25}^{\alpha,0,1} (\alpha\in
    \mathcal{A}), A_{4,25}^{\alpha,\sigma_\xi(1-\alpha),1} (\alpha\in
    \mathcal{A}), \\
    & A_{4,25}^{\alpha,\sigma_\zeta(1-\alpha),\zeta} (\alpha\in
    \mathcal{A}), A_{4,25}^{\alpha,0,\xi} (\alpha\in
    \mathcal{A})
  \end{align*}
  is the list of 4-dimensional nilpotent associative algebras
  (up to isomorphism) if $|F|>3$. If $|F|=3$ then the list is obtained from
  the above one by adding $A_{4,22}$ and erasing
  $A_{4,25}^{\alpha,\sigma_\zeta(1-\alpha),\zeta}$ (which is just one algebra in this
  case). 
\end{prop}

\begin{proof}
The bulk of the list follows directly from the list in Section \ref{sec:list}.
We consider the cases where there is something to do.

Consider the algebra $A_{3,1}$ with basis $a,b,c$. Then $A_{4,8}^{\alpha,\beta}$
is the 1-dimensional central extension of $A_{3,1}$ corresponding to the
cocycle $\theta_{\alpha,\beta}= \Delta_{a,a}+\alpha\Delta_{b,b}+\beta\Delta_{c,c}$.
If we replace $a,b,c$ by $a'= \delta_1 a$, $b'=\delta_2 b$, $c'=\delta_3 c$
then the coefficients of $\Delta_{a,a}$, $\Delta_{b,b}$, $\Delta_{c,c}$ are
multiplied by $\delta_1^2$, $\delta_2^2$, $\delta_3^2$ respectively.
Furthermore, multiplying $\theta_{\alpha,\beta}$ by a nonzero scalar leads to
an isomorphic algebra. By these operations we can reduce to considering
two cocycles: $\theta_{1,1}$, $\theta_{1,\eta}$.
By \cite{omeara}, 62:1 quadratic forms over $F$ are universal, so
there exist $a_{31},a_{32}\in F$ with $a_{31}^2+a_{32}^2 =
\eta$. Now set $a' = a_{32}b+a_{31}c$, $b'=-a_{31}b+a_{32}c$, $c'=\eta a$.
Note that they are linearly independent as $a_{31}^2+a_{32}^2=\eta\neq 0$.
Then we get that $\theta_{1,1} = \eta \Delta_{a',a'} + \eta\Delta_{b',b'}
+\eta^2 \Delta_{c',c'}$. We see that after dividing by $\eta$ we
get $\theta_{1,\eta}$. Hence of the algebras $A_{4,6}^{\alpha,\beta}$ there 
remains only one: $A_{4,6}^{1,1}$.

As noted in Section \ref{sec:A31}, $A_{4,9}^{\alpha,\beta} \cong
A_{4,9}^{\gamma,\delta}$ if and only if $\beta=\delta$ and there are $x,y\in F$
with $x^2+(4\beta-1)y^2 = 4 \tfrac{\alpha}{\gamma}$. If $\beta\neq \tfrac{1}{4}$
then this equation has a solution because quadratic forms are universal.
However, if $\beta=\tfrac{1}{4}$ then the algebras are isomorphic only if
$\alpha$ is a square times $\gamma$.

Consider the algebras $A_{4,25}^{\alpha,\beta,\gamma}$ and set $\sigma = \tfrac{\beta}
{1-\alpha}$. These split into two classes: the first has $\gamma$ equal to
a fixed square (for example 1), the second has $\gamma$ equal to a fixed
non-square (for example $\xi$ or $\zeta$). Each class
again splits into two: the first has $\sigma^2+\gamma$ equal to a fixed
square (for example 1), and the second has $\sigma^2+\gamma$ equal to a
fixed non-square. This leads to the listed algebras.

If $|F|=3$ then the algebra $A_{4,22}$ is added, but
$A_{4,25}^{\alpha,\sigma_\zeta(1-\alpha),\zeta}$ is erased because there are no
non-squares $\zeta\in F$ such that $1-\zeta$ is a square.
\end{proof}

\begin{cor}
  The number of isomorphism classes of 4-dimensional nilpotent associative
  algebras over $F$ is $5q+20$, where $q=|F|$.
\end{cor}  

\subsection{Finite fields of even characteristic}

In this section $F$ is a finite field of even characteristic.

\begin{lem}\label{lem:Hsig}
Let $\sigma \in F$ and set
$$H_\sigma = \{ \tfrac{\sigma uv +v^2}{u^2+v^2} \mid u,v \in F, u\neq v \}.$$
Then $H_\sigma$ is an additive subgroup of $F$. Its index in $F$ is 1 if
$\sigma = 0,1$, and it is 2 otherwise.
\end{lem}

\begin{proof}
It is straightforward to see that $H_\sigma =F$ if $\sigma=0,1$. So suppose
that $\sigma\neq 0,1$. Define $\mathcal{G} = \{ (u,v)\in F\times F \mid
u\neq v\}$. Then $\mathcal{G}$ is a group with group operation $(u_1,v_1)
+(u_2,v_2) = (u_1u_2+v_1v_2,u_1v_2+v_1u_2)$. (Indeed, the neutral element is
$(1,0)$, the inverse of $(u,0)$ is $(u^{-1},0)$, the inverse of $(u,v)$ is
$(u',v')$ with $v' = (u^2v^{-1}+v)^{-1}$, $u'=uv^{-1}v'$ if $v\neq 0$.)
Furthermore, $\tau : \mathcal{G} \to H_\sigma$ given by $\tau(u,v) =
\tfrac{\sigma uv +v^2}{u^2+v^2}$ is a surjective group homomorphism.
Its kernel is
$\{(u,0) \mid u\neq 0\} \cup \{ (u,\sigma u) \mid u\neq 0\}$, which has
$2(q-1)$ elements ($q=|F|$). Now $|\mathcal{G}| = q^2-q$, so $H_\sigma$, being
the image of $\tau$, has $\tfrac{1}{2} q$ elements.
\end{proof}  

\begin{prop}
  Let $\gamma_0$ lie outside the additive subgroup
  $\{t^2+t \mid t \in F\}$ of $F$. For $\sigma \in F\setminus{0,1}$, fix
  $\eta_\sigma$ outside $H_\sigma$. Then
  \begin{align*}
    & A_{4,1}, A_{4,2}, A_{4,3}^1, A_{4,4}^\delta (\delta\in F),
    A_{4,5}, A_{4,6}, A_{4,8}^{1,1}, A_{4,9}^{1,\beta} (\beta\in F), A_{4,10}^\alpha
    (\alpha=0,\gamma_0), A_{4,11},\\
    & A_{4,12}, A_{4,13}, A_{4,14}, A_{4,15}, A_{4,16}, A_{4,17}^\delta (\delta = 0,
    \gamma_0) A_{4,18}^1, A_{4,19}^\delta (\delta=0,\gamma_0), A_{4,20},\\ 
    & A_{4,21}^\delta (\delta\in F, \delta\neq 1), A_{4,26}^{0,0,1}, A_{4,26}^{0,1,1},
    A_{4,26}^{0,\sigma,1} (\sigma \in F\setminus \{0,1\}),
    A_{4,26}^{1+\eta_\sigma^{-1},\sigma \eta_\sigma^{-1},1} (\sigma \in F\setminus \{0,1\}).
  \end{align*}
  is the list of 4-dimensional nilpotent associative algebras over $F$
  (up to isomorphism).
\end{prop}

\begin{proof}
  We only have to comment on the enumeration of the algebras
  $A_{4,26}^{\alpha,\beta,\gamma}$. Obviously we may assume that $\gamma=1$.
  Then $A_{4,26}^{\alpha,\beta,1}$ is isomorphic to $A_{4,26}^{\alpha',\beta',1}$
  if and only if $\sigma=\sigma'$ (where $\sigma = \tfrac{\beta}{1+\alpha}$,
  $\sigma' = \tfrac{\beta'}{1+\alpha'}$) and $\tfrac{1}{1+\alpha} =
  \tfrac{1}{1+\alpha'} +h$ for some $h\in H_\sigma$ (notation as in
  Lemma \ref{lem:Hsig}). By that lemma, for $\sigma=0,1$ we have only one
  algebra with $\tfrac{1}{1+\alpha}=1$ and $\beta=\sigma$. For the other values
  of $\sigma$ we get two algebras, one with $\tfrac{1}{1+\alpha}=1$ and
  $\beta=\sigma$, and one with $\tfrac{1}{1+\alpha}=\eta_\sigma$ and $\beta=
  \sigma(1+\alpha)$.
\end{proof}  

\begin{cor}
  The number of isomorphism classes of 4-dimensional nilpotent associative
  algebras over $F$ is $5q+17$, where $q=|F|$.
\end{cor}

\subsection{The classification over $\R$}

\begin{prop}
  Let $F=\R$. Then
  \begin{align*}
    & A_{4,1}, A_{4,2}, A_{4,3}^\delta (\delta=\pm 1), A_{4,4}^\delta (\delta\in F),
    A_{4,5}, A_{4,6}, A_{4,7}, A_{4,8}^{1,\beta} (\beta=\pm 1), A_{4,9}^{1,\beta}
    (\beta\in F), \\
    & A_{4,9}^{-1,\beta} (\beta\in F, \beta\geq \tfrac{1}{4}),
    A_{4,10}^0, A_{4,11}, A_{4,12}, A_{4,13}, A_{4,14}, A_{4,15}, A_{4,16},
    A_{4,18}^\delta (\delta = 0,\pm 1), \\
    & A_{4,20}
    A_{4,21}^\delta (\delta\in F, \delta\neq -1),
    A_{4,23}^\delta (\delta = \pm 1), 
    A_{4,24}, A_{4,25}^{\alpha,0,\gamma} (\alpha\in
    (-1,1), \gamma=\pm 1),\\
    & A_{4,25}^{\alpha,\sqrt{2}(1-\alpha),-1} (\alpha\in
    (-1,1)). 
  \end{align*}
  is the list of 4-dimensional nilpotent associative algebras
  (up to isomorphism) over $F$.
\end{prop}

\begin{proof}
  For the enumeration of $A_{4,8}^{\alpha,\beta}$, use \cite{jac2}, Chaper V,
  Section 9.
  As remarked in Section \ref{sec:A31}, $A_{4,9}^{\alpha,\beta} \cong
A_{4,9}^{\gamma,\delta}$ if and only if $\beta=\delta$ and there are $x,y\in F$
with $x^2+(4\beta-1)y^2 = 4 \tfrac{\alpha}{\gamma}$. If $\beta < \tfrac{1}{4}$
then this equation always has a solution. So in this case we have one algebra,
$A_{4,9}^{1,\beta}$. If $\beta > \tfrac{1}{4}$, then the equation has a solution
if and only if $\alpha$, $\gamma$ have the same sign. So we obtain two algebras,
$A_{4,9}^{1,\beta}$, $A_{4,9}^{-1,\beta}$.

For the algebras $A_{4,25}$ we remark that we may assume that $\gamma = \pm 1$
and $\sigma^2+\gamma = \pm 1$. If $\sigma^2+\gamma=\gamma$ then $\sigma=0$,
implying $\beta=0$. We have that $A_{4,25}^{\alpha,0,\gamma}\cong A_{4,25}^{\alpha',0,
\gamma}$ if and only if $\alpha=\alpha'$, or $\alpha=(\alpha')^{-1}$. It follows
that by restricting $\alpha$ to the interval $(-1,1)$ we obtain the list of
non-isomorphic algebras $A_{4,25}^{\alpha,0,\gamma}$. If $\gamma=-1$ and
$\sigma^2+\gamma=1$ then we may assume that $\sigma = \sqrt{2}$. Then
again we restrict $\alpha$ to the interval $(-1,1)$, and have
$\beta= \sqrt{2}(1-\alpha)$.
\end{proof}  

\subsection{Algebraically closed fields}

Over algebraically closed fields the enumeration of the algebras is
straightforward (and we leave it to the reader). We remark that if the
characteristic is not 2, then the algebras $A_{4,25}^{\alpha,\beta,\gamma}$
are enumerated as $A_{4,25}^{\alpha,0,1}$, where $\alpha$ runs through a
maximal subset $\mathcal{A}$ of $F$ not containing 1, $-1$ and such that for
$\tau \in \mathcal{A}$, we do not have $\tau^{-1}\in \mathcal{A}$. If $F=\C$,
then for $\mathcal{A}$ we can take the unit circle with the part of the
boundary lying in the upper half plane included.

If the characteristic is 2, then the group $H_{\alpha,\beta,\gamma}$ is all of
$F$, and we may assume that $\alpha=\beta=0$, $\gamma=1$. So here the
class $A_{4,26}$ reduces to one algebra, $A_{4,26}^{0,0,1}$.

\subsection{The classification over $\Q$}

Over $\Q$ we are not able to obtain a very explicit classification.
However, we are able to solve the isomorphism problem. In most cases
it is enough to decide whether a given rational number is a square.
Deciding whether $A_{4,8}^{\alpha,\beta}\cong A_{4,8}^{\gamma,\delta}$ is,
by Lemma \ref{lem:quad}, equivalent to deciding whether the quaternion algebras 
$\left( \tfrac{-\alpha,-\beta}{\Q} \right)$,
$\left( \tfrac{-\gamma,-\delta}{\Q} \right)$ are isomorphic. The latter
question can be decided by computing the sets of places of ramification of
the quaternion algebras (see \cite{vigneras}, Theorem 3.1). Deciding
whether $A_{4,9}^{\alpha,\beta}\cong A_{4,9}^{\gamma,\delta}$ boils down to checking
whether $\beta=\delta$ and whether the curve $x^2+(4\beta-1)y^2 =
4\tfrac{\alpha}{\gamma}$ has a point over $\Q$; the latter can be done using
the methods of \cite{simon}.

\section{The isomorphism problem}\label{sec:isompb}

Given a nilpotent associative algebra of dimension 4, it is possible to follow
the steps in the proof of the classification to obtain the element of
the list of Section \ref{sec:list} to which the given algebra is isomorphic.
We illustrate this in an example. Let
$$A = \langle a,b,c,d \mid a^2=c, b^2=d\rangle$$
(so $A$ is isomorphic to the direct sum of two copies of $A_{2,2}$). Then
$A$ is a 2-dimensional central extension of $A_{2,1}$, so we are in the
situation of Section \ref{sec:noncom2d}. Using the notation in that section
we have $\theta_1 = \left( \begin{smallmatrix}
  1 & 0 \\ 0 & 0 \end{smallmatrix}\right)$,
$\theta_2 = \left(\begin{smallmatrix} 0 & 0 \\ 0 & 1\end{smallmatrix}\right)$.
So by \eqref{eqcoord} $A$ corresponds to the point $a_1=[0,0,1,0,0,0]$.
We first conjugate this to a point with first coordinate nonzero.  
According to the proof of Lemma \ref{lem:A21.1} we can choose $u=v=y=1$, $x=0$
and see that $a_1$ is conjugate to $a_2=[1,1,1,0,0,0]$. If the characteristic
is 2 then we are done (Lemma \ref{lem:A21.2}), and conclude that $A$ is
isomorphic to $A_{4,17}^0$. If the characteristic is not 2 then according to the
proof of Lemma \ref{lem:A21.2} we can choose $u=y=1$, $v=0$, $x=-\tfrac{1}{2}$
and get that $a_2$ is conjugate to $a_3=[1,1,0,0,-\tfrac{1}{4},-\tfrac{1}{4}]$.
By Lemma \ref{lem:A21.6}, $a_3$ is conjugate to $[1,1,0,0,-1,-1]$, which
corresponds to the algebra $A_{4,23}^{-1}$.

\end{document}